\documentclass[twocolumn]{autart} 
\usepackage{graphicx}
\usepackage[T1]{fontenc}
\usepackage{amsmath}
\usepackage{amssymb}
\usepackage{amsfonts}
\usepackage[dvipsnames]{xcolor}
\usepackage{float}
\usepackage{subfigure}
\usepackage{url}
\usepackage{calligra}
\usepackage{enumitem}
\newtheorem{definition}{Definition}

\newtheorem{theorem}{Theorem}

\usepackage{epstopdf}

\begin{document}

\begin{frontmatter}

\title{Event-Triggered Newton Extremum Seeking \\ for Multivariable Optimization}

\thanks[footnoteinfo]{Corresponding author: V.~H.~P.~Rodrigues.}

\author[First]{Victor Hugo Pereira Rodrigues,}
\author[First]{Tiago Roux Oliveira}, ~~~~~~~~~~~~~~~~~~~~~~~~~~~~~~
\author[Second]{Miroslav Krsti{\' c}},  
\author[Third]{Paulo Tabuada} 

\address[First]{Department of Electronics and Telecommunication Engineering,\\ State University of Rio de Janeiro (UERJ), Rio de Janeiro--RJ, Brazil. \\ \mbox{(e-mail: victor.rodrigues@uerj.br, tiagoroux@uerj.br)}}

\address[Second]{Department of Mechanical and Aerospace  Engineering,\\ University of California at San Diego (UCSD), La Jolla--CA, USA. \\ \mbox{(e-mail:  mkrstic@ucsd.edu)}}

\address[Third]{Department of Electrical and Computer Engineering,\\ University of California at Los Angeles (UCLA), Los Angeles--CA, USA. \\ \mbox{(e-mail:  tabuada@ee.ucla.edu)}}

\begin{keyword}                           
Newton-based extremum seeking; Gradient-based extremum seeking; Event-Triggered Control; Discontinuous averaging theory; Multivariable maps; Lyapunov stability.  
\end{keyword}

\begin{abstract}
This paper presents a static event-triggered control strategy for multivariable Newton-based extremum seeking. The proposed method integrates event-triggered actuation into the Newton-based optimization framework to reduce control updates while maintaining rapid convergence to the extremum. Unlike traditional gradient-based extremum seeking, where the convergence rate depends on the unknown Hessian of the cost function, the proposed approach employs a dynamic estimator of the Hessian inverse, formulated as a Riccati equation, enabling user-assignable convergence rates. The event-triggering mechanism is designed to minimize unnecessary actuation updates while preserving stability and performance. Using averaging theory, we establish local stability results and exponential convergence to a neighborhood of the unknown extremum point. Additionally, numerical simulations illustrate the benefits of the proposed approach over gradient-based and continuously actuated Newton-based extremum seeking, showing improved convergence rates and reduced control update frequency, leading to more efficient implementation in real-time optimization scenarios.
\end{abstract}

\end{frontmatter}

\section{Introduction}
In 1922, the concept of extremum seeking was introduced by the French engineer Maurice Leblanc to maintain the maximum efficiency of power transfer from an electrical transmission line to a tram car \cite{L:1922}. Almost eighty years later, a rigorous  stability analysis of extremum seeking feedback was proposed in \cite{KW:2000}, using averaging and singular perturbation methods. This approach optimizes the performance of an unknown dynamic system by employing probing and demodulation periodic signals in order to estimate the gradient and search for the extremum point of the corresponding output cost function.

After that, numerous studies on the topic have been published over the time, including theoretical advances and engineering practice. \cite{AK:2003,Atta2015,Krstic2000,zhang2012extremum}.  In the field of applications, for instance, extremum seeking has been used in autonomous vehicles and mobile sensors \cite{Scheinker2014,SD:2009} as well as optimization of PID controllers applied to Neuromuscular Electrical Stimulation (NMES) \cite{TRO:2019,TRO:2020}. Recent works have also employed extremum seeking to provide policies for Nash equilibrium seeking in non-cooperative games \cite{FKB:2012}, including cases where the players actions are subject to distinct dynamics governed by diffusive Partial Differential Equations (PDEs) \cite{TRoux:2021may} and with delays \cite{TRoux:2021dec,ORKB:2021,ZF:2022,ZFO:2022}. As illustrated in the book \cite{TRoux:2022}, this strategy could be extended to other classes of PDEs. 

On the other hand, due to the rapid development of network technologies, communication between devices and plants have become fast and reliable, reducing wiring costs and simplifying maintenance \cite{ZHGDDYP:2020}. However, one significant disadvantage remains: network bandwidth is limited, and network traffic congestion is unavoidable, leading to delays and packet loss \cite{HNX:2007,Mazenc2022}. To mitigate this problem, event-triggered Control (ETC) can be employed \cite{HJT:2012,Wang2012}. In contrast to periodic sample-data controllers, ETC execute control tasks only when a predefined condition based on the plant’s state is triggered, thereby reducing the consumption of computational resources with guaranteed asymptotic stability properties \cite{Espitia2020,T:2007,Zhang2017}.

In this context, combining real-time optimization by extremum seeking with event-triggered control maybe interesting in order to inherit the advantages of both methods. Oriented for scalar static maps, an event-triggered control scheme for scalar extremum seeking was proposed in \cite{VHPR:2022} and \cite{VHPR:2023b}, where the closed-loop stability of the closed-loop system was guaranteed precluding Zeno behavior. Later on, both static and dynamic triggering scenarios were developed in \cite{VHPR:2023a}, showing the benefits from this combination for the multivariable case as well.

The references cited above \cite{VHPR:2022,VHPR:2023b,VHPR:2023a} employed gradient-based extremum seeking to maximize or minimize the objective functions, where the algorithm's convergence rate speed is proportional to the unknown Hessian (second derivative) of the nonlinear map to be optimized. In 2010, for a single-input case, an estimate of the second derivatives was employed in a Newton-like method \cite{MMB:2010} to reduce sensitivity to the curvature of the plant map. The Newton method requires knowledge of the inverse of the Hessian, which is not a trivial task to provide in model-free optimization. To mitigate this, reference \cite{GKN:2012} proposed a dynamic estimator in the form of a Riccati differential equation for the estimation of the Hessian's inverse matrix, eliminating the dependence on the unknown second derivative in the convergence rate, and making it user-assignable. The key difference between the gradient-based and Newton-based methods is that the former depends on the Hessian, while the latter is independent of it. This offers a significant advantage to the extremum seeking implementation, where the Hessian is unknown.

In this paper, we focus on the design and analysis of a multivariable Newton-based extremum seeking feedback for static maps within an event-triggered framework, aiming to leverage the benefits of improved convergence rates, independent of the Hessian, by employing a Riccati filter to estimate its inverse. We examine the impact of control signal updates and the resulting input-output responses. To address this, a Lyapunov-based  criterion and averaging theory for discontinuous systems are utilized to characterize the stability of the closed-loop system. A numerical example illustrates the theoretical results. 

\textbf{Notation.} Throughout the manuscript, the 2-norm (Euclidean) of vectors and induced norm of matrices are denoted by double bars $\|\cdot\|$ while absolute value of scalar variables are denoted by single bars $|\cdot|$. The terms $\lambda_{\min}(\cdot)$ and $\lambda_{\max}(\cdot)$ denote the minimum and maximum eigenvalues of a given positive definite matrix, respectively. Consider $\varepsilon \in \lbrack -\varepsilon_{0}\,, \varepsilon_{0} \rbrack \subset \mathbb{R}$ and the mappings $\delta_{1}(\varepsilon)$ and $\delta_{2}(\varepsilon)$, where $\delta_{1}: \lbrack -\varepsilon_{0}\,, \varepsilon_{0} \rbrack \to \mathbb{R}$ and $\delta_{2}: \lbrack -\varepsilon_{0}\,, \varepsilon_{0} \rbrack \to \mathbb{R}$. The function $\delta_1(\varepsilon)$ has magnitude of order $\delta_2(\varepsilon)$, denoted by $\delta_{1}(\varepsilon) = \mathcal{O}(\delta_{2}(\varepsilon))$, if there exist positive constants $k$ and $c$ such that $|\delta_{1}(\varepsilon)| \leq k |\delta_{2}(\varepsilon)|$, for all $|\varepsilon|<c$.

\section{Review of the Event-Triggered Multivariable Gradient-based Extremum Seeking}

We consider the following  nonlinear static map
\begin{align}
y(t)&=Q(\theta(t))= Q^{\ast}+\frac{1}{2}(\theta(t)-\theta^{\ast})^{\top}H^{\ast}(\theta(t)-\theta^{\ast}) \label{eq:y_v2} \\
&=Q^{\ast}+\frac{1}{2}\sum_{j=1}^{N}\sum_{k=1}^{N}H_{jk}^{\ast}(\theta_{j}(t)-\theta^{\ast}_{j})(\theta_{k}(t)-\theta^{\ast}_{k})\,, \label{eq:y_v3}
\end{align}
where $Q^{\ast}\in \mathbb{R}$ is the extremum point, $H^{\ast}=H^{\ast \top} \in \mathbb{R}^{n \times n}$ is the Hessian matrix, $\theta^{\ast} \in \mathbb{R}^{n}$ is the optimizer, $\theta(t)\in \mathbb{R}^{n}$ is the input map. The cost function $Q(\theta(t))$ is the map to be optimized and its parameters in (\ref{eq:y_v2}) are not explicitly known; however, we have access to measurements of $y(t) \in \mathbb{R}$ and can adjust $\theta(t)$. 

If the Hessian matrix is positive definite (i.e., $H^{\ast} > 0$), the map (\ref{eq:y_v2}) is convex and attains a minimum at $\theta = \theta^{\ast}$. Conversely, if $H^{\ast} < 0$, the function becomes concave and the extremum at $\theta^{\ast}$ corresponds to a maximum. Therefore, by evaluating the sign of $H^{\ast}$, one can directly infer whether the extremum seeking algorithm should be interpreted as solving a minimization or maximization problem.

Although (\ref{eq:y_v2}) is a polynomial in $\theta$, it cannot be identified from a finite number of input/output pairs $(\theta, y) $ due to the real-time nature of extremum seeking optimization. Furthermore, extremum seeking is a versatile approach that can handle any unknown analytic function $Q(\theta)$, which can be expressed as a Taylor series expansion around a point $\theta^{\ast}$, where the function attains a minimum or maximum. This assumption enables a local quadratic approximation of the nonlinear map $ Q(\theta)$, serving as a foundational principle and justifying our focus on quadratic maps. Therefore, while the analysis was conducted for a quadratic map, our approach is not limited to quadratic functions $Q(\theta)$.

\subsection{Classical Gradient-based Extremum Seeking}

\begin{figure*}[h!]
\centering
\includegraphics[width=13cm]{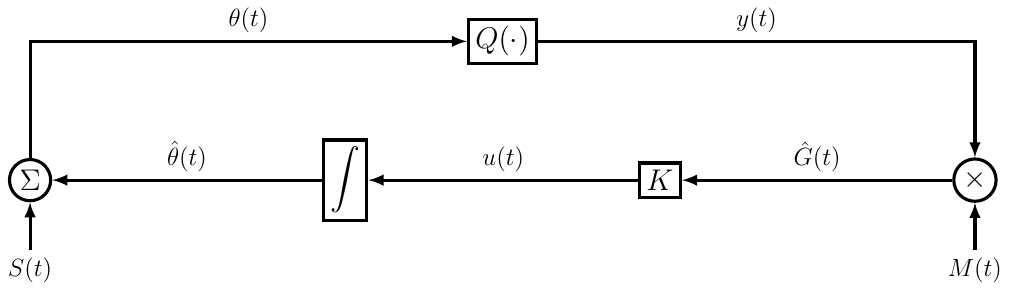}
\caption{Multivariable Gradient-based Extremum Seeking.}
\label{fig:BD_GradientES_v2}
\end{figure*}

The multivariable gradient-based Extremum Seeking (GradientES) approach for this multivariable static map is illustrated in Fig.~\ref{fig:BD_GradientES_v2}. In this scheme, the feedback or adaptation gain is
\begin{align}
K=\text{diag}\left\{K_{1}\,,K_{2}\,,\ldots\,,K_{n}\right\}\,, \label{eq:K}
\end{align}
while the dither signals and demodulations are defined as (see \cite{GKN:2012,K:2014}:
\begin{align}
S(t)&= \left[a_{1}\sin\left(\omega_1 t\right),\ldots,a_{i}\sin\left(\omega_i t\right),\ldots,a_{n}\sin\left(\omega_n t\right)\right]^{\top}\!\!, \label{eq:S_v1} \\
M(t)&=2\left[\frac{\sin\left(\omega_1 t\right)}{a_{1}},\ldots,\frac{\sin\left(\omega_i t\right)}{a_{i}},\ldots,\frac{\sin\left(\omega_n t\right)}{a_{n}}\right]^{\top}, \label{eq:M_v1}
\end{align}
with nonzero amplitudes $a_{i}$. Moreover, the probing frequencies $\omega_{i}$'s can be selected as
\begin{align}
\omega_{i}=\omega_{i}'\omega \,, \quad i \in \left\{1,\ldots\,,n\right\}\,, \label{eq:omegai_event}
\end{align}
where $\omega$ is a positive constant and $\omega_{i}'$ is a rational number. In \cite{GKN:2012}, the following assumption is originally introduced.
\begin{assum}\label{assumption_w}
The probing frequencies satisfy
\begin{align}
\omega'_{i} 	\notin \left\{\omega'_{j}\,,~\frac{1}{2}(\omega'_{j}+\omega'_{k})\,,~\omega'_{j}+2\omega'_{k}\,,~\omega'_{k}\pm\omega'_{l}\right\}\,, \label{eq:omega_iNotIn}
\end{align}
for all $i$, $j$, $k$ and $l$.
\end{assum} 

The scheme in Fig.~\ref{fig:BD_GradientES_v2} highlights the core components of the ES architecture: the gradient estimation mechanism driven by the periodic perturbation-demodulation technique and the adaptation block that adjusts the input vector using the estimated gradient in order to estimate $\theta^{\ast}$.

Moreover, if Assumption~\ref{assumption_w} holds, an accurate estimate of the $i$-th component of the unknown gradient vector can be obtained by means of $\hat{G}(t) = M(t)y(t) \in \mathbb{R}^{n}$, whose $i$-th component can be expressed as
\begin{align}
\hat{G}_{i}(t)&=\frac{2}{a_{i}}\sin(\omega_{i}t)y(t)\,. \label{eq:Gi}
\end{align}
As mentioned above, the output of the integrator, the vector $\hat{\theta}(t)\in \mathbb{R}^{n}$, give us an estimate of $\theta^{*} \in \mathbb{R}^{n}$ such that the \textit{estimation error} is defined by: 
\begin{align}
\tilde{\theta}(t)&:=\hat{\theta}(t)-\theta^{*}\,.\label{eq:tildeThetai_v1}
\end{align}
Still referring to Fig.~\ref{fig:BD_GradientES_v2}, the input to the nonlinear map $Q(\theta)$ in (\ref{eq:y_v2}) is given by $\theta(t) = \hat{\theta}(t) + S(t)$. For analysis purposes, and using (\ref{eq:S_v1}) and (\ref{eq:tildeThetai_v1}), the $i$-th component of the input vector $\theta(t)$ can be expressed in terms of the $i$-th component of the unknown estimation error vector $\tilde{\theta}(t)$, and the $i$-th component of the unknown optimal parameter vector $\theta^{\ast}$, as well as the $i$-th component of the known dither signal $S(t)$, as follows
\begin{align}
\theta_{i}(t)&=\tilde{\theta}_{i}(t)+a_{i}\sin\left(\omega_i t\right)+\theta_{i}^{*}\,.\label{eq:thetai_v1}
\end{align}
Therefore, by plugging (\ref{eq:y_v3}),  (\ref{eq:M_v1}) and (\ref{eq:thetai_v1}) into (\ref{eq:Gi}), the $i$-th component of the estimate of the unknown gradient vector can be expressed as
\begin{align}
&\hat{G}_{i}(t)=\frac{1}{a_{i}}\sin(\omega_{i}t)\tilde{\theta}^{\top}(t)H^{\ast}\tilde{\theta}(t)\nonumber \\
&\qquad~~~~~+\sum_{j=1}^{n}\mbox{\calligra H}_{~~ij}^{~~~\ast}(t)\tilde{\theta}_{j}(t) +\Delta_{i}(t) \,, \label{eq:hatG_20250522_1} \\
&\mbox{\calligra H}_{~~ij}^{~~~\ast}(t):=H_{ij}^{\ast}+\Delta\mbox{\calligra H}_{~~ij}^{~~~\ast}(t)\,, \label{eq:calligraH_ij_20250522_1} \\
&\Delta\mbox{\calligra H}_{~~ij}^{~~~\ast}(t):=-H_{ij}^{\ast}\cos(2\omega_{i}t)\nonumber \\
&\qquad~~~~~~~~~~~~~+\sum_{\substack{k=1 \\ k\neq i}}^{n}H_{kj}^{\ast}\frac{a_{k}}{a_{i}}\cos((\omega_{i}-\omega_{k})t) \nonumber \\
&\qquad~~~~~~~~~~~~~-\sum_{\substack{k=1 \\ k\neq i}}^{n}H_{kj}^{\ast}\frac{a_{k}}{a_{i}}\cos((\omega_{i}+\omega_{k})t)\,, \label{eq:DeltaCalligraH_ij_20250522_1} \\
%
&\Delta_{i}(t):=\frac{2Q^{\ast}}{a_{i}}\sin(\omega_{i}t) \nonumber \\
&\qquad~~~~~~+\frac{1}{4}\sum_{j=1}^{n}\sum_{k=1}^{n}H_{jk}^{\ast}\frac{a_{j}a_{k}}{a_{i}}\sin((\omega_{i}+\omega_{j}-\omega_{k})t) \nonumber \\
&\qquad~~~~~~-\frac{1}{4}\sum_{j=1}^{n}\sum_{k=1}^{n}H_{jk}^{\ast}\frac{a_{j}a_{k}}{a_{i}}\sin((\omega_{i}-\omega_{j}+\omega_{k})t) \nonumber \\
&\qquad~~~~~~-\frac{1}{4}\sum_{j=1}^{n}\sum_{k=1}^{n}H_{jk}^{\ast}\frac{a_{j}a_{k}}{a_{i}}\sin((\omega_{i}+\omega_{j}+\omega_{k})t) \nonumber \\
%
&\qquad~~~~~~-\frac{1}{4}\sum_{j=1}^{n}\sum_{k=1}^{n}H_{jk}^{\ast}\frac{a_{j}a_{k}}{a_{i}}\sin((\omega_{i}-\omega_{j}-\omega_{k})t)\,. \label{eq:Delta_i_20250222_1}
\end{align}
Hence, by using (\ref{eq:calligraH_ij_20250522_1})--(\ref{eq:Delta_i_20250222_1}), we can define the time-varying matrices $\mbox{\calligra H}^{~~~\ast}(t) \in \mathbb{R}^{n \times n}$, $\Delta \! \mbox{\calligra H}^{~~~\ast}(t) \in \mathbb{R}^{n \times n}$, and the time-varying vector $\Delta(t) \in \mathbb{R}^{n }$, as follows
\begin{align}
\mbox{\calligra H}^{~~~\ast}(t) &:= H^{\ast}+\Delta \! \mbox{\calligra H}^{~~~\ast}(t)\,, \label{eq:calligraH}\\
\Delta \! \mbox{\calligra H}^{~~~\ast}(t)&:= \begin{bmatrix}
													\Delta \! \mbox{\calligra H}^{~~~\ast}_{~~~11}(t) & \Delta \! \mbox{\calligra H}^{~~~\ast}_{~~~12}(t) & \ldots & \Delta \! \mbox{\calligra H}^{~~~\ast}_{~~~1n}(t) \\
													\Delta \! \mbox{\calligra H}^{~~~\ast}_{~~~21}(t) & \Delta \! \mbox{\calligra H}^{~~~\ast}_{~~~22}(t) & \ldots & \Delta \! \mbox{\calligra H}^{~~~\ast}_{~~~2n}(t) \\
													\vdots                         & \vdots                         & \ddots & \vdots                         \\
													\Delta \! \mbox{\calligra H}^{~~~\ast}_{~~~n1}(t) & \Delta \! \mbox{\calligra H}^{~~~\ast}_{~~~n2}(t) & \ldots & \Delta \! \mbox{\calligra H}^{~~~\ast}_{~~~nn}(t) \\
												 \end{bmatrix} \,, \label{eq:DeltaCalligraH} \\
						 \Delta(t) &:= \begin{bmatrix}
													\Delta_{1}(t) \,, 
													\Delta_{2}(t) \,,
													\ldots\,,
													\Delta_{n}(t)
												 \end{bmatrix}^{\top} \,. \label{eq:Delta}
\end{align}
Then, by using (\ref{eq:M_v1}) and (\ref{eq:calligraH})--(\ref{eq:Delta}), we can express (\ref{eq:hatG_20250522_1}) in the compact form
\begin{align}
\hat{G}(t)&\mathbb{=}M(t)\frac{1}{2}\tilde{\theta}^{\top}(t)H^{\ast}\tilde{\theta}(t)\mathbb{+}\left(H^{\ast}\mathbb{+}\Delta \! \mbox{\calligra H}^{~~~\ast}\!\!(t)\right)\tilde{\theta}(t)\mathbb{+}\Delta(t), \label{eq:hatG_20240302_2}
\end{align}
where $\Delta \! \mbox{\calligra H}~(t)$, defined in (\ref{eq:DeltaCalligraH}), and $\Delta(t)$, given in (\ref{eq:Delta}), are time-varying matrix and vector, respectively, both with zero mean.

Since the term $\tilde{\theta}^{\top}(t)H^{\ast}\tilde{\theta}(t)$ is quadratic in $\tilde{\theta}(t)$ and, therefore, may be neglected in a local analysis \cite{AK:2003},  the gradient estimate (\ref{eq:hatG_20240302_2}) can simply be rewritten as 
\begin{align}
\hat{G}(t)&=\left(H^{\ast}+\Delta \! \mbox{\calligra H}^{~~~\ast}(t)\right)\tilde{\theta}(t)+\Delta(t)\,. \label{eq:hatG_20240302_3}
\end{align}

On the other hand, from the time-derivative of (\ref{eq:tildeThetai_v1}) and the classical GradientES scheme depicted in Fig.~\ref{fig:BD_GradientES_v2}, the dynamics that governs $\hat{\theta}(t)$, as well as $\tilde{\theta}(t)$, is given by
\begin{align}
\frac{d\tilde{\theta}}{dt}(t)&=\frac{d\hat{\theta}}{dt}(t)=u(t) \label{eq:dtildeThetadt_20250206_1}\,, 
\end{align}
with $u(t)= [u_{1}(t)\,,u_{2}(t)\,,\ldots\,,u_{n}(t)]^{\top} \in \mathbb{R}^{n}$. Moreover, by using (\ref{eq:dtildeThetadt_20250206_1}), the time-derivative of (\ref{eq:hatG_20240302_3}) is given by 
\begin{align}
\frac{d\hat{G}}{dt}(t)&\mathbb{=}\left(H^{\ast}\mathbb{+}\Delta \! \mbox{\calligra H}^{~~~\ast}\!(t)\right)u(t)\mathbb{+}\frac{d\Delta \! \mbox{\calligra H}^{~~~\ast}}{dt}\!(t)\tilde{\theta}(t)\mathbb{+}\frac{d \Delta}{dt}(t), \label{eq:dhatGdt_20250206_1}
\end{align}
where, from (\ref{eq:calligraH}) and (\ref{eq:Delta}), it is easy to verify that $\dfrac{d\Delta \! \mbox{\calligra H}^{~~~\ast}}{dt}(t)$ as well as $\dfrac{d \Delta}{dt}(t)$ have null mean values as well.

\subsection{Static Gradient-Feedback Tuning Law}

We assume the static gradient-feedback control law, being, for all $t\geq 0$, 
\begin{align}
u(t)=K\hat{G}(t) \,, \quad \forall t\geq 0 \label{eq:u}
\end{align}
such that $KH^{\ast}$ is Hurwitz. 

We assume that the control law (\ref{eq:u}) stabilizes the system at the corresponding average equilibrium $\hat{G}_{\rm{av}} \equiv 0$ with local exponential convergence. This means that any trajectory of $\hat{\theta}(t)$ starting sufficiently near to the extremum point $\theta^{\ast}$ will converge to a neighborhood of it exponentially, with uniform decay and overshoot bounds. It simply formalizes the idea that the control law is designed for local stabilization, regardless of the specific value of $\theta$. The control law (\ref{eq:u}) does not rely on detailed knowledge of the mapping (\ref{eq:y_v3}) or the system dynamics (\ref{eq:dtildeThetadt_20250206_1}) and (\ref{eq:dhatGdt_20250206_1}). Its design is model-independent in that sense, which is particularly advantageous in scenarios where precise models are unavailable or hard to obtain.

\subsection{Event-Triggered Multivariable Gradient-Based Extremum Seeking} \label{sec:ET-GradientES}

\begin{figure*}[h!]
\centering
\includegraphics[width=17cm]{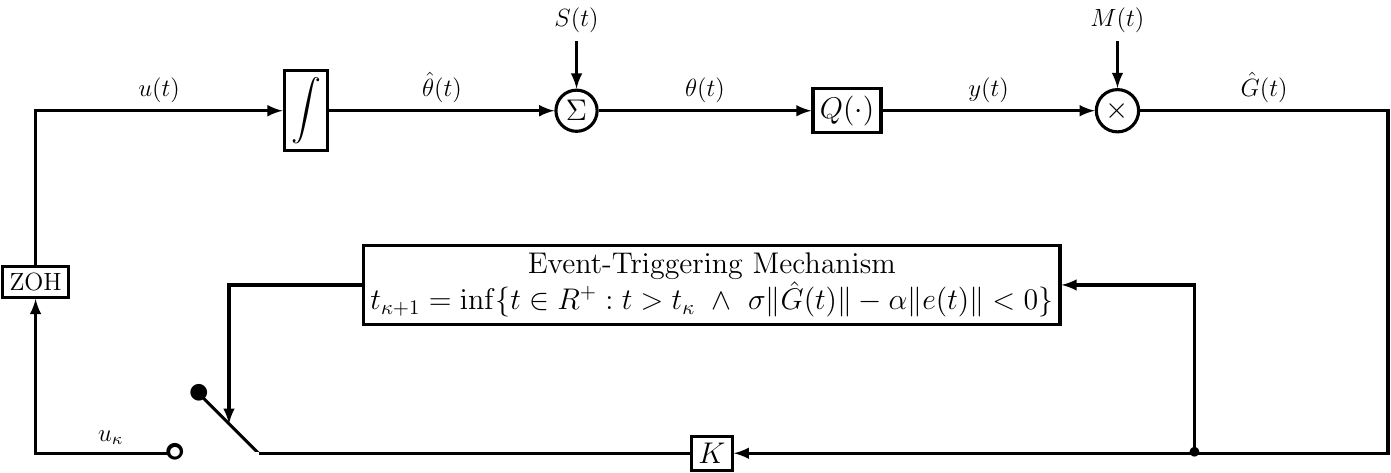}
\caption{Event-Triggered Multivariable Gradient-based Extremum Seeking.}
\label{fig:BD_ETGradientES}
\end{figure*}

In our previous work~\cite{VHPR:2023a}, a static event-triggered implementation was proposed to emulate the behavior of the continuous-time multivariable gradient-based extremum seeking controller (\ref{eq:u}), while significantly reducing the frequency of control updates. As shown in Fig.~\ref{fig:BD_ETGradientES}, the key idea lies in monitoring a suitably designed error-dependent condition to determine when the control law needs to be updated, rather than enforcing updates at every time instant.

To formally characterize the stability properties of the resulting event-triggered control loop, the following theorem was established in \cite{VHPR:2023a}:

\begin{theorem} \label{thm:ETGradientES}
\noindent
Consider the dynamics of the gradient estimate given by equation~(\ref{eq:dhatGdt_20250206_1}) with control law (\ref{eq:u}) emulated by $u(t)=K\hat{G}(t_{k})$, for all $t \in [t_k, t_{k+1})$, where $t_k$ denotes the most recent update time. The next update instant $t_{k+1}$ is determined by the event-triggering mechanism: $t_{0}=0$ and $t_{k+1} = \inf\left\{ t > t_k \; \middle| \; \sigma \|\hat{G}(t)\|- \alpha  \|e(t)\| < 0 \right\}$, along with the error vector defined as $e(t) := \hat{G}(t_k) - \hat{G}(t)$, where $\sigma \in (0,1)$ and $\alpha>0$ are design parameters. Assume that the matrix product $H^{\ast}K$ is Hurwitz such there exists the symmetric positive definite matrices $P = P^\top > 0$ and $Q = Q^\top > 0$ satisfying the Lyapunov equation $K^\top H^{\ast \top} P + P H^{\ast} K = -Q$, and let Assumptions~\ref{assumption_w} hold. Then, for the defined period $T = 2\pi \times \mathrm{LCM}\left\{ \frac{1}{\omega_i} \right\}$, frequency $\omega=\frac{2\pi}{T}$ sufficiently large, and $\alpha > \frac{2 \| P H^{\ast} K \|}{\lambda_{\min}(Q)}$, one has
\begin{small}
   \begin{align}
\|\theta(t)\mathbb{-}\theta^{\ast}\|&\mathbb{\leq} \sqrt{\frac{\lambda_{\max}(H^{\ast \top}PH^{\ast})}{\lambda_{\min}(H^{\ast \top}PH^{\ast})}} \exp\left(\mathbb{-}\frac{\lambda_{\min}(Q)}{2\lambda_{\max}(P)}(1\mathbb{-}\sigma)t\right)\nonumber \\
&\quad\mathbb{\times}\|\theta(0) - \theta^{\ast}\|+\mathcal{O}\left(a+\frac{1}{\omega}\right)\,, \label{eq:normTheta_thm2} \\ 
|y(t)\mathbb{-}Q^{\ast}|&\mathbb{\leq} \|H^{\ast}\|\frac{\lambda_{\max}(H^{\ast \top}PH^{\ast})}{\lambda_{\min}(H^{\ast \top}PH^{\ast})}\exp\left(\!\!\mathbb{-}\frac{\lambda_{\min}(Q)}{2\lambda_{\max}(P)}(1\mathbb{-}\sigma)t\right)  \nonumber \\
&\quad \mathbb{\times}\|\theta(0) \mathbb{-} \theta^{\ast}\|^{2}+\mathcal{O}\left(  a^{2} +\frac{1}{\omega^{2}}\right)\,, \label{eq:normY_thm2}
\end{align} 
\end{small}
where $a=\sqrt{\sum_{i=1}^{n}a_{i}^{2}}$, with $a_i$ defined in (\ref{eq:S_v1}) and (\ref{eq:M_v1}), and the constants $M_{\theta}$, and $M_{y}$ depending on the initial condition $\theta(0)$. In addition, there exists a lower bound  $\tau^{\ast}$ for the inter-execution interval $t_{k+1}-t_{k}$  for all $k \in \mathbb{N}$ precluding the Zeno behavior.
\end{theorem}

The closest proof of Theorem~\ref{thm:ETGradientES} is developed in \cite{VHPR:2023a}.

Despite the advantages of reducing control updates and preserving stability through event-triggered execution, the multivariable gradient-based extremum seeking framework suffers from several fundamental limitation that motivate the design of a multivariable Newton-based alternative. A key drawback lies in its sensitivity to the unknown curvature of the cost function: the convergence rate of gradient-based methods is inherently tied to the eigenstructure of the Hessian matrix, which is typically unknown and may lead to slow convergence along poorly conditioned directions. This is especially problematic in multivariable settings, where highly elongated level sets cause the parameter trajectories to follow inefficient, curved ``steepest descent'' paths. 

In order to overcome this limitation, the next section introduces a new event-triggered multivariable Newton-based extremum seeking approach, which leverages a dynamic estimation of the inverse Hessian matrix to decouple convergence performance from unknown curvature information, enabling faster, more direct convergence with fewer control updates.

\section{Event-Triggered Multivariable Newton-based Extremum Seeking} \label{sec:eventTiggred}

\begin{figure*}[h!]
\centering
\includegraphics[width=12cm]{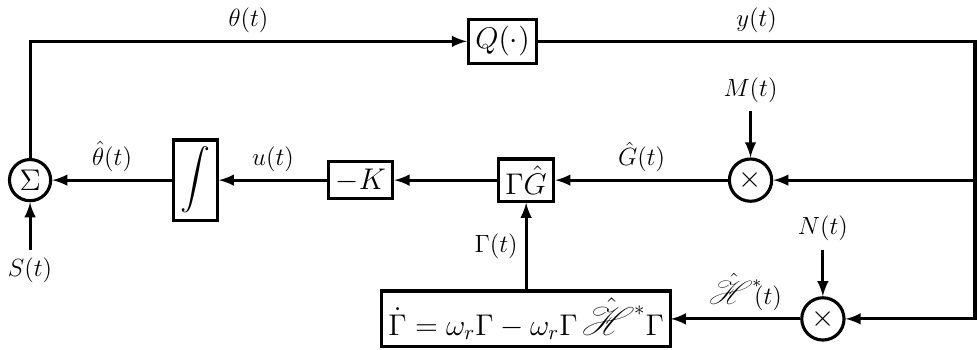}
\caption{Multivariable Newton-based Extremum Seeking.}
\label{fig:BD_NewtonES}
\end{figure*}

Fig.~\ref{fig:BD_NewtonES} depicts the Newton-based extremum seeking algorithm applied to a static multi-input map. By comparing the block diagrams shown in Fig.~\ref{fig:BD_GradientES_v2} and Fig.~\ref{fig:BD_NewtonES}, it becomes evident that in the Newton-based extremum seeking control, both the gradient estimation and the estimation error of the optimal point are computed identically to those in the gradient-based extremum seeking control scheme. Therefore, by following steps given from (\ref{eq:Gi})--(\ref{eq:dhatGdt_20250206_1}), in the Newton-based extremum seeking strategy, we have
\begin{align}
\tilde{\theta}(t)&\mathbb{=}\hat{\theta}(t)-\theta^{*}\,,\label{eq:tildeThetai_20251121_v1} \\
\frac{d\tilde{\theta}}{dt}(t)&\mathbb{=}\frac{d\hat{\theta}}{dt}(t)=u(t) \label{eq:dtildeThetadt_20251121_1}\,, \\
\hat{G}(t)&\mathbb{=}\left(H^{\ast}+\Delta \! \mbox{\calligra H}^{~~~\ast}(t)\right)\tilde{\theta}(t)+\Delta(t)\,, \label{eq:hatG_20251121_3} \\
\frac{d\hat{G}}{dt}(t)&\mathbb{=}\left(H^{\ast}\mathbb{+}\Delta \! \mbox{\calligra H}^{~~~\ast}\!(t)\right)u(t)\mathbb{+}\frac{d\Delta \! \mbox{\calligra H}^{~~~\ast}}{dt}\!(t)\tilde{\theta}(t)\mathbb{+}\frac{d \Delta}{dt}(t), \label{eq:dhatGdt_20251127_1}
\end{align}
where, from (\ref{eq:DeltaCalligraH_ij_20250522_1}), (\ref{eq:Delta_i_20250222_1}), (\ref{eq:DeltaCalligraH}) and (\ref{eq:Delta}), it is easy to verify that $\Delta \! \mbox{\calligra H}^{~~~\ast}(t)$, $\dfrac{d\Delta \! \mbox{\calligra H}^{~~~\ast}}{dt}(t)$ as well as $\Delta(t)$ and $\dfrac{d \Delta}{dt}(t)$ have null mean values.

Moreover, the Riccati filter---driven by a positive real frequency $\omega_r$---represents the key additional component in the Newton-based extremum seeking architecture depicted in Fig.~\ref{fig:BD_NewtonES}. This filter is responsible for dynamically estimating the inverse Hessian and distinguishes the Newton-based approach from its gradient-based counterpart. The overall algorithm relies on two essential mechanisms: gradient estimation and curvature compensation via the inverse Hessian estimate.

\subsection{Continuous-time Newton-based Extremum Seeking: Hessian Estimate and its Inverse}

The first is the \emph{perturbation structure}, defined by the time-varying matrix $N(t)\in \mathbb{R}^{n \times n}$, with elements 
\begin{align}
    N_{ii}(t)&\mathbb{:=}\mathbb{-}\frac{8}{a_{i}^2}\cos(2\omega_{i}t)\,, \label{eq:Nii}\\
    N_{ij}(t)&\mathbb{:=}\frac{2}{a_{i}a_{j}}\cos((\omega_{i}\mathbb{-}\omega_{j})t)\mathbb{-}\frac{2}{a_{i}a_{j}}\cos((\omega_{i}\mathbb{+}\omega_{j})t),~ j\neq i , \label{eq:Nij}
\end{align}
which introduce carefully chosen input oscillations that allow for the estimation of second-order derivative information of the unknown map (\ref{eq:y_v2})--namely, the Hessian matrix $H^{\ast}$, by means of
\begin{align}
    \widehat{\!\!\!\!\!\!\!\!\mbox{\calligra H}}^{~~~\ast}(t) &:= N(t)y(t)\,. \label{eq:calligraHHat_20250527}
\end{align}
Therefore, by plugging (\ref{eq:y_v3}),  (\ref{eq:Nii}) and (\ref{eq:Nij}) into (\ref{eq:calligraHHat_20250527}), the components of the estimate  ~~~$\widehat{\!\!\!\!\!\!\!\!\mbox{\calligra H}}^{~~\ast}_{~~ii}(t)$ and ~~~$\widehat{\!\!\!\!\!\!\!\!\mbox{\calligra H}}^{~~\ast}_{~~ij}(t)$ of the unknown Hessian matrix $H^{\ast}$ are given by
\begin{align}
    \widehat{\!\!\!\!\!\!\!\!\mbox{\calligra H}}^{~~\ast}_{~~ii}(t)&\mathbb{=} H^{\ast}_{ii} \mathbb{+}\Delta~~~~\! \widehat{\mbox{\calligra \!\!\!\!\!\!\!\! H}}^{~~\ast}_{~~ii}(t) \nonumber \\
    &\mathbb{+}N_{ii}(t)\left[\frac{1}{2}\tilde{\theta}^{\top}(t)H^{\ast}\tilde{\theta}(t)\mathbb{+}S^{\top}(t)H^{\ast}\tilde{\theta}(t)\right]\,, \label{eq:calligraHiiHat} \\
    \Delta~~~~\! \widehat{\mbox{\calligra \!\!\!\!\!\!\!\! H}}^{~~\ast}_{~~ii}(t)&= -H^{\ast}_{ii}\cos(2\omega_{i}t)+N_{ii}(t)Q^{\ast} \nonumber  \\
    & \!\!\!\!\!\!\!\!\!\!\!\!\!\!\!\!\!\!\!\!\!\!\!\!\!\! -2\sum_{\substack{k=1 \\ k\neq i}}^{n}H^{\ast}_{kk}\frac{a_{k}^{2}}{a_{i}^{2}}\cos(2\omega_{i}t) \nonumber \\
    & \!\!\!\!\!\!\!\!\!\!\!\!\!\!\!\!\!\!\!\!\!\!\!\!\!\! +\sum_{\substack{k=1 \\ k\neq i}}^{n}H^{\ast}_{kk}\frac{a_{k}^{2}}{a_{i}^{2}}\left[\cos(2(\omega_{i}-\omega_{k})t)+\cos(2(\omega_{i}+\omega_{k})t)\right] \nonumber \\
    &\!\!\!\!\!\!\!\!\!\!\!\!\!\!\!\!\!\!\!\!\!\!\!\!\!\! -2\sum_{k=1}^{n}\sum_{l=k+1}^{n}H^{\ast}_{kl}\frac{a_{k}a_{l}}{a_{i}^2}\cos((2\omega_{i}+\omega_{k}-\omega_{l})t) \nonumber \\
        &\!\!\!\!\!\!\!\!\!\!\!\!\!\!\!\!\!\!\!\!\!\!\!\!\!\! -2\sum_{k=1}^{n}\sum_{l=k+1}^{n}H^{\ast}_{kl}\frac{a_{k}a_{l}}{a_{i}^2}\cos((2\omega_{i}-\omega_{k}+\omega_{l})t) \nonumber  \\
         &\!\!\!\!\!\!\!\!\!\!\!\!\!\!\!\!\!\!\!\!\!\!\!\!\!\! +2\sum_{k=1}^{n}\sum_{l=k+1}^{n}H^{\ast}_{kl}\frac{a_{k}a_{l}}{a_{i}^2}\cos((2\omega_{i}+\omega_{k}+\omega_{l})t) \nonumber \\
    &\!\!\!\!\!\!\!\!\!\!\!\!\!\!\!\!\!\!\!\!\!\!\!\!\!\! +2\sum_{k=1}^{n}\sum_{l=k+1}^{n}H^{\ast}_{kl}\frac{a_{k}a_{l}}{a_{i}^2}\cos((2\omega_{i}-\omega_{k}-\omega_{l})t)\,, \label{eq:DeltaCalligraHiiHat} \\
    \widehat{\!\!\!\!\!\!\!\!\mbox{\calligra H}}^{~~\ast}_{~~ij}(t)&\mathbb{=} H^{\ast}_{ij} \mathbb{+}\Delta~~~~\! \widehat{\mbox{\calligra \!\!\!\!\!\!\!\! H}}^{~~\ast}_{~~ij}(t) \nonumber \\
    &\mathbb{+}N_{ij}(t)\left[\frac{1}{2}\tilde{\theta}^{\top}(t)H^{\ast}\tilde{\theta}(t)\mathbb{+}S^{\top}(t)H^{\ast}\tilde{\theta}(t)\right]\,, \label{eq:calligraHijHat} 
\end{align}
\begin{align}
    \Delta~~~~\! \widehat{\mbox{\calligra \!\!\!\!\!\!\!\! H}}^{~~\ast}_{~~ij}(t)&= -H_{ij}^{\ast}\cos(2\omega_{i}t)-H_{ij}^{\ast}\cos(2\omega_{j}t)\nonumber \\
    &\!\!\!\!\!\!\!\!\!\!\!\!\!\!\!\!\!\!\!\!\!\!\!\!\!\! +\frac{H_{ij}^{\ast}}{2}\cos(2(\omega_{i}-\omega_{j})t)+\frac{H_{ij}^{\ast}}{2}\cos(2(\omega_{i}+\omega_{j})t) \nonumber \\
    &\!\!\!\!\!\!\!\!\!\!\!\!\!\!\!\!\!\!\!\!\!\!\!\!\!\! +\frac{1}{2}\sum_{\substack{k=1 \\ k\neq i}}^{n}H^{\ast}_{kk}\frac{a_{k}^{2}}{a_{i}a_{j}}\left[\cos(2(\omega_{i}-\omega_{j})t)+\cos(2(\omega_{i}+\omega_{j})t)\right] \nonumber  \\
    &\!\!\!\!\!\!\!\!\!\!\!\!\!\!\!\!\!\!\!\!\!\!\!\!\!\! -\frac{1}{2}\sum_{\substack{k=1 \\ k\neq i}}^{n}H^{\ast}_{kk}\frac{a_{k}^{2}}{a_{i}a_{j}}\cos((\omega_{i}-\omega_{j}+2\omega_{k})t) \nonumber  \\
    &\!\!\!\!\!\!\!\!\!\!\!\!\!\!\!\!\!\!\!\!\!\!\!\!\!\! -\frac{1}{2}\sum_{\substack{k=1 \\ k\neq i}}^{n}H^{\ast}_{kk}\frac{a_{k}^{2}}{a_{i}a_{j}}\cos((\omega_{i}-\omega_{j}-2\omega_{k})t) \nonumber  \\
    &\!\!\!\!\!\!\!\!\!\!\!\!\!\!\!\!\!\!\!\!\!\!\!\!\!\! +\frac{1}{2}\sum_{\substack{k=1 \\ k\neq i}}^{n}H^{\ast}_{kk}\frac{a_{k}^{2}}{a_{i}a_{j}}\cos((\omega_{i}+\omega_{j}+2\omega_{k})t) \nonumber  \\
    &\!\!\!\!\!\!\!\!\!\!\!\!\!\!\!\!\!\!\!\!\!\!\!\!\!\! +\frac{1}{2}\sum_{\substack{k=1 \\ k\neq i}}^{n}H^{\ast}_{kk}\frac{a_{k}^{2}}{a_{i}a_{j}}\cos((\omega_{i}+\omega_{j}-2\omega_{k})t) \nonumber  \\
    &\!\!\!\!\!\!\!\!\!\!\!\!\!\!\!\!\!\!\!\!\!\!\!\!\!\! +\frac{1}{2}\sum_{\substack{k=1 \\ k\neq i}}^{n}\sum_{l=k+1}^{n}H^{\ast}_{kl}\frac{a_{k}a_{l}}{a_{i}a_{j}}\cos((\omega_{i}-\omega_{j}+\omega_{k}-\omega_{l})t) \nonumber \\
    &\!\!\!\!\!\!\!\!\!\!\!\!\!\!\!\!\!\!\!\!\!\!\!\!\!\! +\frac{1}{2}\sum_{\substack{k=1 \\ k\neq i}}^{n}\sum_{l=k+1}^{n}H^{\ast}_{kl}\frac{a_{k}a_{l}}{a_{i}a_{j}}\cos((\omega_{i}-\omega_{j}-\omega_{k}+\omega_{l})t) \nonumber \\
    &\!\!\!\!\!\!\!\!\!\!\!\!\!\!\!\!\!\!\!\!\!\!\!\!\!\! -\frac{1}{2}\sum_{\substack{k=1 \\ k\neq i}}^{n}\sum_{l=k+1}^{n}H^{\ast}_{kl}\frac{a_{k}a_{l}}{a_{i}a_{j}}\cos((\omega_{i}-\omega_{j}+\omega_{k}+\omega_{l})t) \nonumber \\
    &\!\!\!\!\!\!\!\!\!\!\!\!\!\!\!\!\!\!\!\!\!\!\!\!\!\! -\frac{1}{2}\sum_{\substack{k=1 \\ k\neq i}}^{n}\sum_{l=k+1}^{n}H^{\ast}_{kl}\frac{a_{k}a_{l}}{a_{i}a_{j}}\cos((\omega_{i}-\omega_{j}-\omega_{k}-\omega_{l})t) \nonumber  \\
        &\!\!\!\!\!\!\!\!\!\!\!\!\!\!\!\!\!\!\!\!\!\!\!\!\!\! -\frac{1}{2}\sum_{\substack{k=1 \\ k\neq i}}^{n}\sum_{l=k+1}^{n}H^{\ast}_{kl}\frac{a_{k}a_{l}}{a_{i}a_{j}}\cos((\omega_{i}+\omega_{j}+\omega_{k}-\omega_{l})t) \nonumber  \\
        &\!\!\!\!\!\!\!\!\!\!\!\!\!\!\!\!\!\!\!\!\!\!\!\!\!\! -\frac{1}{2}\sum_{\substack{k=1 \\ k\neq i}}^{n}\sum_{l=k+1}^{n}H^{\ast}_{kl}\frac{a_{k}a_{l}}{a_{i}a_{j}}\cos((\omega_{i}+\omega_{j}-\omega_{k}+\omega_{l})t) \nonumber \\
             &\!\!\!\!\!\!\!\!\!\!\!\!\!\!\!\!\!\!\!\!\!\!\!\!\!\! +\frac{1}{2}\sum_{\substack{k=1 \\ k\neq i}}^{n}\sum_{l=k+1}^{n}H^{\ast}_{kl}\frac{a_{k}a_{l}}{a_{i}a_{j}}\cos((\omega_{i}+\omega_{j}+\omega_{k}+\omega_{l})t) \nonumber \\
             &\!\!\!\!\!\!\!\!\!\!\!\!\!\!\!\!\!\!\!\!\!\!\!\!\!\! +\frac{1}{2}\sum_{\substack{k=1 \\ k\neq i}}^{n}\sum_{l=k+1}^{n}H^{\ast}_{kl}\frac{a_{k}a_{l}}{a_{i}a_{j}}\cos((\omega_{i}+\omega_{j}-\omega_{k}-\omega_{l})t) \nonumber  \\
              &\!\!\!\!\!\!\!\!\!\!\!\!\!\!\!\!\!\!\!\!\!\!\!\!\!\! +\frac{1}{2}\sum_{\substack{k=i+1 \\ k\neq j}}^{n}H^{\ast}_{ik}\frac{a_{k}}{a_{j}}\cos((2\omega_{i}-\omega_{j}-\omega_{k})t) \nonumber \\ 
        &\!\!\!\!\!\!\!\!\!\!\!\!\!\!\!\!\!\!\!\!\!\!\!\!\!\! +\frac{1}{2}\sum_{\substack{k=i+1 \\ k\neq j}}^{n}H^{\ast}_{ik}\frac{a_{k}}{a_{j}}\cos((-\omega_{j}+\omega_{k})t) \nonumber
    \end{align}
    \begin{align} 
        &\!\!\!\!\!\!\!\!\!\!\!\!\!\!\!\!\!\!\!\!\!\!\!\!\!\! -\frac{1}{2}\sum_{\substack{k=i+1 \\ k\neq j}}^{n}H^{\ast}_{ik}\frac{a_{k}}{a_{j}}\cos((2\omega_{i}+\omega_{j}+\omega_{k})t) \nonumber \\
    &\!\!\!\!\!\!\!\!\!\!\!\!\!\!\!\!\!\!\!\!\!\!\!\!\!\! -\frac{1}{2}\sum_{\substack{k=i+1 \\ k\neq j}}^{n}H^{\ast}_{ik}\frac{a_{k}}{a_{j}}\cos((\omega_{j}+\omega_{k})t)+N_{ij}(t)Q^{\ast}\,. \label{eq:DeltaCalligraHijHat} 
\end{align}
Thus, by using (\ref{eq:calligraHiiHat})--(\ref{eq:DeltaCalligraHijHat}), we can define the time-varying matrix $\Delta~~~~\! \widehat{\mbox{\calligra \!\!\!\!\!\!\!\! H}}^{~~\ast}(t) \in \mathbb{R}^{n \times n}$ such that
\begin{align}
\Delta~~~~\! \widehat{\mbox{\calligra \!\!\!\!\!\!\!\! H}}^{~~~\ast}(t)&\mathbb{:=} \begin{bmatrix}
													\Delta~~~~\! \widehat{\mbox{\calligra \!\!\!\!\!\!\!\! H}}^{~~~\ast}_{~~~11}(t) & \Delta~~~~\! \widehat{\mbox{\calligra \!\!\!\!\!\!\!\! H}}^{~~~\ast}_{~~~12}(t) & \ldots & \Delta~~~~\! \widehat{\mbox{\calligra \!\!\!\!\!\!\!\! H}}^{~~~\ast}_{~~~1n}(t) \\
													\Delta~~~~\! \widehat{\mbox{\calligra \!\!\!\!\!\!\!\! H}}^{~~~\ast}_{~~~21}(t) & \Delta~~~~\! \widehat{\mbox{\calligra \!\!\!\!\!\!\!\! H}}^{~~~\ast}_{~~~22}(t) & \ldots & \Delta~~~~\! \widehat{\mbox{\calligra \!\!\!\!\!\!\!\! H}}^{~~~\ast}_{~~~2n}(t) \\
													\vdots                         & \vdots                         &    \ddots    & \vdots                         \\
													\Delta~~~~\! \widehat{\mbox{\calligra \!\!\!\!\!\!\!\! H}}^{~~~\ast}_{~~~n1}(t) & \Delta~~~~\! \widehat{\mbox{\calligra \!\!\!\!\!\!\!\! H}}^{~~~\ast}_{~~~n2}(t) & \ldots & \Delta~~~~\! \widehat{\mbox{\calligra \!\!\!\!\!\!\!\! H}}^{~~~\ast}_{~~~nn}(t) \\
												 \end{bmatrix}, \label{eq:DeltaCalligraHhat} 
\end{align}
and (\ref{eq:calligraHHat_20250527}) can be written as

\begin{align}
   \widehat{\!\!\!\!\!\!\!\!\mbox{\calligra H}}^{~~~\ast}(t) &:= H^{\ast}+\Delta~~~~\! \widehat{\mbox{\calligra \!\!\!\!\!\!\!\! H}}^{~~~\ast}(t) \nonumber \\
   &\quad+N(t)\left[\frac{1}{2}\tilde{\theta}^{\top}(t)H^{\ast}\tilde{\theta}(t)\mathbb{+}S^{\top}(t)H^{\ast}\tilde{\theta}(t)\right]\,. \label{eq:calligraHhat_20250528_v1}
\end{align}
Since the term $\tilde{\theta}^{\top}(t)H^{\ast}\tilde{\theta}(t)$ is quadratic in $\tilde{\theta}(t)$ and, therefore, may be neglected in a local analysis \cite{AK:2003},  the Hessian estimate (\ref{eq:calligraHhat_20250528_v1}) can be rewritten as 
\begin{align}
   ~~~~\widehat{\!\!\!\!\!\!\!\!\mbox{\calligra H}}^{~~~\ast}(t) &:= H^{\ast}\mathbb{+}\Delta~~~~\! \widehat{\mbox{\calligra \!\!\!\!\!\!\!\! H}}^{~~~\ast}(t) \mathbb{+}S^{\top}(t)H^{\ast}\tilde{\theta}(t)N(t) . \label{eq:calligraHhat_20250528_v2}
\end{align}

The second component is a \emph{Riccati-based adaptation law}, 
\begin{align}
    \frac{d\Gamma}{dt}(t):=\omega_{r}\Gamma(t)-\omega_{r}\Gamma(t)~~~~\widehat{\!\!\!\!\!\!\!\!\mbox{\calligra H}}^{~~~\ast}\!\!(t)\Gamma(t)\,, \label{eq:dotGamma}
\end{align}
designed to compute an estimate of the \emph{inverse of the Hessian}. A key feature of this Riccati-based approach is its robustness: it ensures that a usable inverse estimate is available even when the Hessian approximation is close to singular zero values or it is poorly conditioned. Together, these components enable the implementation of a Newton-like optimization scheme in real time, without requiring explicit knowledge of the cost function’s derivatives.

The Riccati equation in (\ref{eq:dotGamma}) admits two equilibrium solutions: $\Gamma^{\ast} = 0_{n \times n}$ and $\Gamma^{\ast} = H^{\ast-1}$, assuming the estimate $~~~~~~\widehat{\!\!\!\!\!\!\!\!\mbox{\calligra H}}^{~~~\ast}\!\!(t)$ converges to a constant. For $\omega_r > 0$, the zero equilibrium is unstable, while linearizing the dynamics around $H^{\ast}$ yields a Jacobian matrix equal to $-\omega_r I_{n}$, indicating that this equilibrium is locally exponentially stable. Consequently, the Riccati dynamics are guaranteed to converge—after a transient phase—to the inverse of the estimated Hessian matrix, provided that~~ $~~\widehat{\!\!\!\!\!\!\!\!\mbox{\calligra H}}^{~~~\ast}\!\!(t)$ reasonably approximates the true Hessian $H^{\ast}$, for more details, see \cite{GKN:2012}. 

While precisely characterizing the region of attraction of the stable equilibrium $\Gamma^{\ast} = H^{\ast-1}$ is generally difficult, a conservative yet practical estimate is inversely proportional to the largest eigenvalue of $H^{\ast}$. This bound provides a useful guideline for understanding how the Hessian estimate affects the convergence behavior of the Riccati-based inverse computation.

\subsection{Continuous-Time Tuning Law}

For the Newton-based algorithm depicted in Fig.~\ref{fig:BD_NewtonES}, the corresponding tuning law is given by
\begin{align}
u(t)=-K\Gamma(t)\hat{G}(t) \,, \quad \forall t\geq 0\,. \label{eq:uNewton}
\end{align}
Thus, in order to analyze the closed-loop system, we use the estimation error (\ref{eq:tildeThetai_v1}), the gradient estimate (\ref{eq:hatG_20240302_3}) and define Hessian estimation error 
\begin{align}
    \tilde{\Gamma}(t):=\Gamma(t)-H^{\ast-1}\,, \label{eq:tildeGamma}
\end{align}
such that the closed-loop system, in error variables, by using (\ref{eq:calligraHHat_20250527})--(\ref{eq:calligraHhat_20250528_v2}), is given by 
\begin{align}
\frac{d\hat{G}}{dt}(t)&= \mathbb{-}\left(H^{\ast}KH^{\ast\mathbb{-}1}\mathbb{+}H^{\ast}K\tilde{\Gamma}(t) \right. \nonumber \\
&\!\!\left.\mathbb{+}\Delta \! \mbox{\calligra H}^{~~~\ast}\!(t)KH^{\ast\mathbb{-}1}\mathbb{+}\Delta \! \mbox{\calligra H}^{~~~\ast}\!(t)K\tilde{\Gamma}(t) \right)\hat{G}(t) \nonumber \\
&\mathbb{+}\frac{d\Delta \! \mbox{\calligra H}^{~~~\ast}}{dt}\!(t)\tilde{\theta}(t)\mathbb{+}\frac{d \Delta}{dt}(t), \label{eq:dotHatG_20250603_v1} \\
\frac{d\tilde{\theta}}{dt}(t)&=\mathbb{-}K\tilde{\theta}(t)\mathbb{-}K\tilde{\Gamma}(t)H^{\ast}\tilde{\theta}(t)\mathbb{-}KH^{\ast\mathbb{-}1}\Delta(t) \nonumber \\
&\mathbb{-}KH^{\ast\mathbb{-}1}\Delta \! \mbox{\calligra H}^{~~~\ast}\!(t)\tilde{\theta}(t)\mathbb{-}K\tilde{\Gamma}(t)\Delta \! \mbox{\calligra H}^{~~~\ast}\!(t)\tilde{\theta}(t)\label{eq:dtildeThetadt_20250529_1}\,,\\
\frac{d\tilde{\Gamma}}{dt}(t)&=\mathbb{-}\omega_{r}\tilde{\Gamma}(t)\mathbb{+}\omega_{r}\Upsilon(t,\tilde{\theta},\tilde{\Gamma})\,, \label{eq:dotTildeGamma_20250529_1} \\
\Upsilon(t,\tilde{\theta},\tilde{\Gamma})&=\mathbb{-}\tilde{\Gamma}(t)H^{\ast}\tilde{\Gamma}(t) \nonumber \\
&\mathbb{-}\tilde{\Gamma}(t)(\Delta~~~~\! \widehat{\mbox{\calligra \!\!\!\!\!\!\!\! H}}^{~~~\ast}(t)\mathbb{+}S^{\top}(t)H^{\ast}\tilde{\theta}(t)N(t))H^{\ast\mathbb{-}1} \nonumber \\
&\mathbb{-}H^{\ast\mathbb{-}1}(\Delta~~~~\! \widehat{\mbox{\calligra \!\!\!\!\!\!\!\! H}}^{~~~\ast}(t)\mathbb{+}S^{\top}(t)H^{\ast}\tilde{\theta}(t)N(t))\tilde{\Gamma}(t) \nonumber \\
&\mathbb{-}\tilde{\Gamma}(t)(\Delta~~~~\! \widehat{\mbox{\calligra \!\!\!\!\!\!\!\! H}}^{~~~\ast}(t)\mathbb{+}S^{\top}(t)H^{\ast}\tilde{\theta}(t)N(t))\tilde{\Gamma}(t) \nonumber \\
&\mathbb{-}H^{\ast\mathbb{-}1}(\Delta~~~~\! \widehat{\mbox{\calligra \!\!\!\!\!\!\!\! H}}^{~~~\ast}(t)\mathbb{+}S^{\top}(t)H^{\ast}\tilde{\theta}(t)N(t))H^{\ast\mathbb{-}1}\,.\label{eq:Upsilon_20251007_1}
\end{align}

Notice that, from (\ref{eq:S_v1}) and (\ref{eq:Nii})--(\ref{eq:calligraHhat_20250528_v2}), the time-varying scalar $S^{\top}(t)H^{\ast}\tilde{\theta}(t)$ as well as the time-varying matrices $N(t)$ and $\Delta~~~ \widehat{\mbox{\calligra \!\!\!\!\!\!\!\! H}}^{~~\ast}(t)$ have null mean values. Moreover,  $\tilde{\Gamma}(t) H^{\ast} \tilde{\theta}(t)$ in (\ref{eq:dtildeThetadt_20250529_1}) introduces a nonlinear coupling between the errors variables, as it is quadratic in the variables $\tilde{\Gamma}(t)$ and $\tilde{\theta}(t)$. Similarly, the term $\tilde{\Gamma}(t) H^{\ast} \tilde{\Gamma}(t)$ in (\ref{eq:Upsilon_20251007_1}) is quadratic in $\tilde{\Gamma}(t)$ , capturing how this error influences itself through second-order interactions in dynamics (\ref{eq:dotTildeGamma_20250529_1}).

Hence, the average and linearized closed-loop system around the equilibrium point being given by
\begin{align}
\frac{d\hat{G}_{\rm{av}}}{dt}(t)&=H^{\ast}(-K)H^{\ast-1}\hat{G}_{\rm{av}}(t)\label{eq:dhatGdt_20250529_2}\,,\\
\frac{d\tilde{\theta}_{\rm{av}}}{dt}(t)&=-K\tilde{\theta}_{\rm{av}}(t)\label{eq:dtildeThetadt_20250529_2}\,,\\
    \frac{d\tilde{\Gamma}_{\rm{av}}}{dt}(t)&=-\omega_{r}\tilde{\Gamma}_{\rm{av}}(t)\,. \label{eq:dotTildeGamma_20250529_2}
\end{align}

From (\ref{eq:dhatGdt_20250529_2}) and (\ref{eq:dtildeThetadt_20250529_2}), the symmetric matrix $H^{\ast}=H^{\ast \top}$ has definite signal, therefore, it full rank and there exists the inverse matrix $H^{\ast -1}$ while $K=K^{\top}$ is a symmetric positive definite diagonal matrix thus although the matrix $H^{\ast}KH^{\ast-1}$ is not necessarily positive definite or even symmetric, it has the same eigenvalues as $K$ since $\det(H^{\ast}KH^{\ast-1}-\lambda I_{n})=\det(H^{\ast}KH^{\ast-1}-\lambda H^{\ast}H^{\ast-1})=\det(H^{\ast}(K-\lambda I_{n})H^{\ast-1})=\det(H^{\ast})\det(K-\lambda I_{n})\det(H^{\ast-1})=\det(H^{\ast})\det(K-\lambda I_{n})\det^{-1}(H^{\ast})=\det(K-\lambda I_{n})$. The state transformation $H^{\ast}KH^{\ast-1}$ represents the same linear operator as $K$, but expressed in a different basis, the basis of $H^{\ast}$. Therefore, the following assumption is additionally considered throughout the paper.

\begin{assum}\label{assumption_lyapunovEq}
The matrices $(-K)$ and $H^{\ast}(-K)H^{\ast-1}$ are Hurwitz such that for any given $Q=Q^{\top}>0$ there exist symmetric positive definite matrices $P_{1}=P_{1}^{\top}>0$ and $P_{2}=P_{2}^{\top}>0$ satisfying the Lyapunov equations
\begin{align}
(-K)^{\top}P_{1}+P_{1}(-K)&=-Q \,, \label{eq:LyapEq_1} \\
(H^{\ast}(-K)H^{\ast-1})^{\top}P_{2}+P_{2}(H^{\ast}(-K)H^{\ast-1})&=-Q \,. \label{eq:LyapEq_2} 
\end{align}
\end{assum}

This property marks a key advantage of the Newton-based extremum seeking algorithm over its gradient-based counterpart. While gradient-based methods have convergence rates inherently tied to the (typically unknown) Hessian matrix $H^{\ast}$, the Newton-based approach allows for the designer to \emph{explicitly assign} the convergence characteristics through the appropriate choice of the gain matrix $K$, in (\ref{eq:dhatGdt_20250529_2}) and (\ref{eq:dtildeThetadt_20250529_2}), and the excitation frequency $\omega_r$ in (\ref{eq:dotTildeGamma_20250529_2}). As a result, tuning these parameters offers direct control over the stability and speed of convergence, independently of the curvature of the cost function.

\subsection{Emulation of the Continuous-Time Tuning Law}

As discussed in Section~\ref{sec:ET-GradientES}, the control updates are executed at a sequence of time instants $(t_{\kappa})_{\kappa \in \mathbb{N}}$, which are not predetermined at uniform intervals, but instead generated by an event-triggering mechanism. This mechanism is designed to ensure the closed-loop system retains both stability and robustness. 

\begin{figure*}[h!]
\centering
\includegraphics[width=17cm]{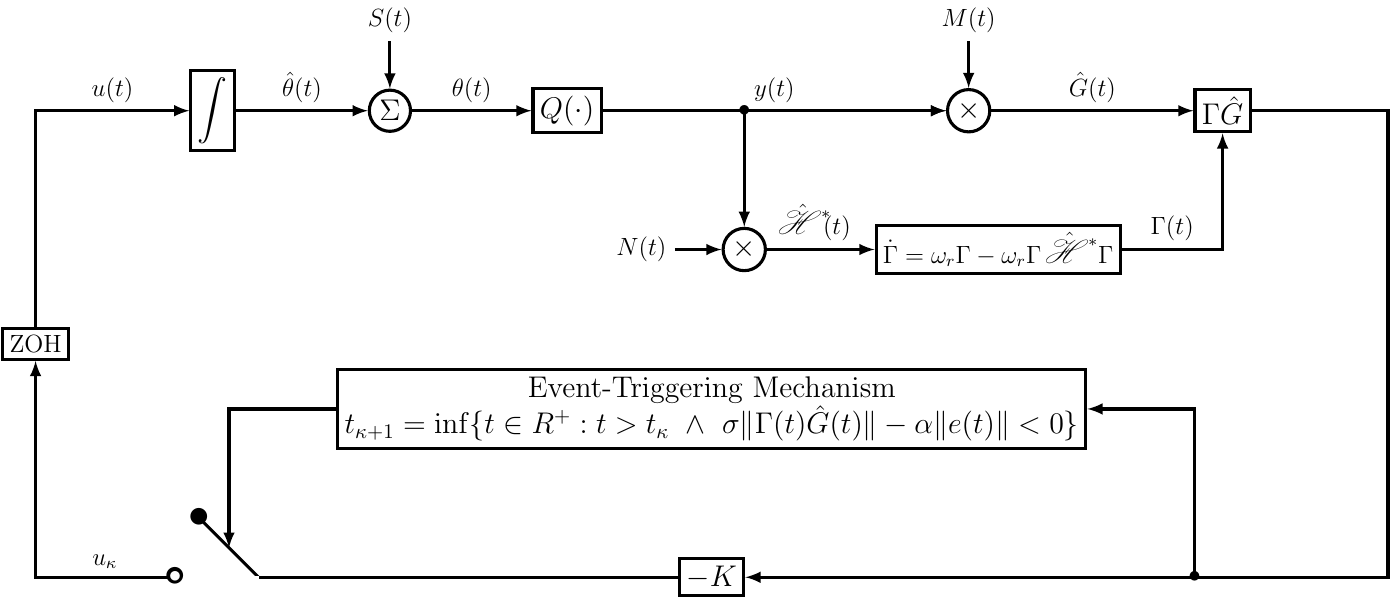}
\caption{Event-Triggered Multivariable Newton-based Extremum Seeking.}
\label{fig:BD_NewtonES_v5}
\end{figure*}

The event-triggering mechanism for multivariable Newton-based extremum seeking operates by continuously monitoring the evolution of the product between the inverse Hessian estimate $\Gamma(t)$ and the gradient estimate $\hat{G}(t)$, and determining whether a control update is required based on a prescribed deviation criterion. In other words, the control input $u(t)$ is only updated at discrete time instants $\{t_\kappa\}_{\kappa \in \mathbb{N}}$, which are generated by a condition that compares the current product $\Gamma(t)\hat{G}(t)$ to its most recently sampled value $\Gamma(t_{\kappa})\hat{G}(t_{\kappa})$. If the difference between these two quantities exceeds a predefined threshold, an event is triggered and the control law is re-evaluated. 

Assuming no delays in the sensor-to-controller and controller-to-actuator communication paths, the control input is held constant between updates using a zero-order hold (ZOH) mechanism. That is, the control input computed at the last triggering instant \(t_\kappa\) is applied continuously over the interval \([t_\kappa, t_{\kappa+1})\), ensuring that the control effort is only refreshed when necessary, thereby reducing unnecessary actuation and improving implementation efficiency such that
\begin{align}
u(t)\!=\!-K\Gamma(t_{\kappa})\hat{G}(t_{\kappa})\,, \quad \forall t \!\in\! \lbrack t_{\kappa}\,, t_{\kappa+1}\phantom{(}\!\!) \,, \quad \kappa\!\in\! \mathbb{N}  \,, \label{eq:u_20250602_v1}
\end{align}
and, since the control updates are made at discrete times rather than continuously, we introduce the error 
\begin{align}
e(t):=\Gamma(t_{\kappa})\hat{G}(t_{\kappa})-\Gamma(t)\hat{G}(t) \,, ~ \forall t \in \lbrack t_{\kappa}\,, t_{\kappa+1}\phantom{(}\!\!) \,,~ \kappa\in \mathbb{N} \,. \label{eq:e_event_20250602_v1}
\end{align}
Therefore, by using (\ref{eq:hatG_20240302_3}), (\ref{eq:tildeGamma}) and (\ref{eq:e_event_20250602_v1}), for all $ t \in \lbrack t_{k}\,, t_{k+1}\phantom{(}\!\!)$, $ k\in \mathbb{N} $, the event-triggered multivariable Newton-based control law (\ref{eq:u_20250602_v1}) can be rewritten as
\begin{align}
u(t)&=-KH^{\ast-1}\hat{G}(t)-Ke(t)-K\tilde{\Gamma}(t)\hat{G}(t)  \,, \label{eq:u_20250602_v2} \\
&=-K\tilde{\theta}(t)-Ke(t)-K\tilde{\Gamma}(t)H^{\ast}\tilde{\theta}(t)-KH^{\ast-1}\Delta(t) \nonumber \\
&\quad-KH^{\ast-1}\Delta \! \mbox{\calligra H}^{~~~\ast}\!(t)\tilde{\theta}(t)-K\tilde{\Gamma}(t)\Delta \! \mbox{\calligra H}^{~~~\ast}\!(t)\tilde{\theta}(t)  \,. \label{eq:u_20250602_v3}
\end{align}
Now, plugging (\ref{eq:u_20250602_v2}) into (\ref{eq:dhatGdt_20251127_1}) and (\ref{eq:u_20250602_v3}) into (\ref{eq:dtildeThetadt_20251121_1}), we arrive at the following Input-to-State Stable (ISS) \cite{K:2002} representations for the dynamics of $\hat{G}(t)$ and $\tilde{\theta}$ with respect to the time-varying error matrix $\tilde{\Gamma}(t)$ in (\ref{eq:tildeGamma}) and the error vector $e(t)$ in (\ref{eq:e_event_20250602_v1}), for all $t\in[t_{\kappa},t_{\kappa+1})$: 
\begin{align}
\frac{d\hat{G}}{dt}(t)&\mathbb{=} \mathbb{-}\left(H^{\ast}KH^{\ast-1}\mathbb{+}H^{\ast}K\tilde{\Gamma}(t)\mathbb{+}\Delta \! \mbox{\calligra H}^{~~~\ast}\!(t)KH^{\ast-1} \right. \nonumber \\
&\quad\!\!\left.\mathbb{+}\Delta \! \mbox{\calligra H}^{~~~\ast}\!(t)K\tilde{\Gamma}(t) \right)\hat{G}(t)\mathbb{+}\frac{d\Delta \! \mbox{\calligra H}^{~~~\ast}}{dt}\!(t)\tilde{\theta}(t) \nonumber \\
&\quad\mathbb{-}\left(\!H^{\ast}\mathbb{+}\Delta \! \mbox{\calligra H}^{~~~\ast}\!(t)\!\right)Ke(t)\mathbb{+}\frac{d \Delta}{dt}(t), \label{eq:dotHatG_20250602_v1} \\
\frac{d\tilde{\theta}}{dt}(t)&\mathbb{=}\mathbb{-}K\left(I_{n}\mathbb{+}\tilde{\Gamma}(t)H^{\ast}\mathbb{+}H^{\ast-1}\Delta \! \mbox{\calligra H}^{~~~\ast}\!(t)\right. \nonumber \\
&\quad\left.\!\!\mathbb{+}\tilde{\Gamma}(t)\Delta \! \mbox{\calligra H}^{~~~\ast}\!(t)\right)\tilde{\theta}(t)-Ke(t)-KH^{\ast-1}\Delta(t)\,, \label{eq:dotTildeTheta_20250602_v1}
\end{align}
where $I_{n} \in \mathbb{R}^{n \times n}$ is an identity matrix.

In the following sections, we present a static event-triggered mechanism for the nonautonomous and nonlinear time-varying closed-loop system given by (\ref{eq:dotTildeGamma_20250529_1}), (\ref{eq:dotHatG_20250602_v1}) and (\ref{eq:dotTildeTheta_20250602_v1}). In this framework, data transmission occurs only when the norm of the measurement error vector (\ref{eq:e_event_20250602_v1}) exceeds a predefined threshold. This triggering mechanism operates without continuous communication, thereby reducing the communication burden while preserving system stability.

\subsection{Event-Triggering Mechanism} 

Definition~\ref{def:staticEvent} illustrates how the small design parameter $\sigma \in (0,1)$ and the error signal $e(t)$—representing the deviation between the product  $\Gamma(t)\hat{G}(t)$ and its last broadcasted value are leveraged to construct a static event-triggered mechanism for multivariable Newton-based extremum seeking. The resulting strategy triggers the recalculation of the control law, as defined in equation~(\ref{eq:u_20250602_v1}), and updates the zero-order hold (ZOH) actuator only when needed, as shown in the block diagram of Fig.~\ref{fig:BD_NewtonES_v5}. By embedding this triggering logic into the Newton-based optimization loop—where convergence speed is decoupled from the unknown Hessian and made user-assignable—the closed-loop system preserves asymptotic stability and guaranteed convergence, as supported by the theoretical results in \cite{HJT:2012}.

\begin{definition}[\small{Event-Triggering Condition}] \label{def:staticEvent}
The  multivariable Newton-based extremum seeking controller with static-event-triggered condition consists of two components:
\begin{enumerate}
	\item  A sequence of increasing times $\mathcal{I}$ such that $\mathcal{I}=\{t_{0}\,, t_{1}\,, t_{2}\,,\ldots\}$ with $t_{0}=0$, generated under the following rules:
		\begin{itemize}[leftmargin=-0.1in]
			\item \!\!\!\! If $\left\{t \in\mathbb{R}^{+}: t \mathbb{>} t_{\kappa} ~ \mathbb{\wedge} ~ \sigma\|\Gamma(t)\hat{G}(t)\| \mathbb{-}\alpha \|e(t)\| \mathbb{<} 0 \right\} \mathbb{=} \emptyset$, then the set of the times of the events is $\mathcal{I}=\{t_{0}\,, t_{1}\,, \ldots, t_{\kappa}\}$.
			\item \!\!\!\! If $\left\{t \in\mathbb{R}^{+}: t \mathbb{>} t_{\kappa} ~ \mathbb{\wedge} ~ \sigma\|\Gamma(t)\hat{G}(t)\| \mathbb{-}\alpha\|e(t)\| \mathbb{<} 0 \right\} \mathbb{\neq} \emptyset$, the next event time, consisting of the static event-triggering mechanism, is given by
				\begin{align}
					\!\!\!\!t_{\kappa+1}&\mathbb{=}\inf\left\{t \in\mathbb{R}^{+}\mathbb{:}~ t \mathbb{>} t_{\kappa} ~ \mathbb{\wedge} ~ \sigma\|\Gamma(t)\hat{G}(t)\| \mathbb{-}\alpha\|e(t)\| \mathbb{<} 0 \right\}\!\!, \label{eq:tk+1_event}
				\end{align} 
        \end{itemize}
	\item A feedback control action (\ref{eq:u_20250602_v1}) updated at time instants $t_{\kappa}$.
\end{enumerate}  
\end{definition}

\subsection{Rescaling of Time} \label{sec:timeScale}

Now, we  introduce a suitable time scale to carry out  the stability analysis of the closed-loop system. From (\ref{eq:omegai_event}), one can notice that the  frequencies in (\ref{eq:S_v1}), (\ref{eq:M_v1}), (\ref{eq:Nii}) and (\ref{eq:Nij}), as well as their combinations, are rational. Furthermore, there exists a time period $T$ such that 
\begin{align}
T&= 2\pi \times \text{LCM}\left\{\frac{1}{\omega_{i}}\right\}\,, \quad \forall i \left\{1\,,2\,,\ldots\,,n\right\}\,, \label{eq:T}
\end{align}
 where LCM denotes the least common multiple such that it is possible to define the time-scale for the dynamics  (\ref{eq:dtildeThetadt_20250206_1}) and (\ref{eq:dhatGdt_20250206_1}) with the transformation $\bar{t}=\omega t$, where 
\begin{align}
\omega&:=\frac{2\pi}{T}\,. \label{eq:omega_event_1}
\end{align}
 Hence, the system (\ref{eq:dotTildeGamma_20250529_1}), (\ref{eq:dotHatG_20250602_v1}) and (\ref{eq:dotTildeTheta_20250602_v1}) can be rewritten as, $\forall t \in \lbrack t_{\kappa}\,, t_{\kappa+1}) \,,  ~\kappa\in \mathbb{N}$, 
\begin{align}
\frac{d\hat{G}}{d\bar{t}}(\bar{t})&\mathbb{=} \mathbb{-}\frac{1}{\omega}\left(H^{\ast}KH^{\ast-1}\mathbb{+}H^{\ast}K\tilde{\Gamma}(\bar{t})\mathbb{+}\Delta \! \mbox{\calligra H}^{~~~\ast}\!(\bar{t})KH^{\ast-1} \right. \nonumber \\
&\quad\!\!\left.\mathbb{+}\Delta \! \mbox{\calligra H}^{~~~\ast}\!(\bar{t})K\tilde{\Gamma}(\bar{t}) \right)\hat{G}(\bar{t})\mathbb{+}\frac{d\Delta \! \mbox{\calligra H}^{~~~\ast}}{d\bar{t}}\!(\bar{t})\tilde{\theta}(\bar{t}) \nonumber \\
&\quad\mathbb{-}\frac{1}{\omega}\left(\!H^{\ast}\mathbb{+}\Delta \! \mbox{\calligra H}^{~~~\ast}\!(\bar{t})\!\right)Ke(\bar{t})\mathbb{+}\frac{d \Delta}{d\bar{t}}(\bar{t}), \label{eq:dhatGdt_20250206_3} \\
\frac{d\tilde{\theta}}{d\bar{t}}(\bar{t})&\mathbb{=}\mathbb{-}\frac{1}{\omega}K\left(\! I_{n}\mathbb{+}\tilde{\Gamma}(\bar{t})H^{\ast}\mathbb{+}H^{\ast-1}\Delta \! \mbox{\calligra H}^{~~~\ast}\!(\bar{t})\right. \nonumber \\
&\quad\left.\!\!\mathbb{+}\tilde{\Gamma}(\bar{t})\Delta \! \mbox{\calligra H}^{~~~\ast}\!(\bar{t})\!\right)\tilde{\theta}(\bar{t})\mathbb{-}\frac{1}{\omega}Ke(\bar{t})\mathbb{-}\frac{1}{\omega}KH^{\ast\mathbb{-}1}\Delta(\bar{t}), \label{eq:dtildeThetadt_20250206_3} \\
\frac{d\tilde{\Gamma}}{d\bar{t}}(\bar{t})&=\!\mathbb{-}\frac{\omega_{r}}{\omega}\tilde{\Gamma}(\bar{t})\mathbb{+}\frac{\omega_{r}}{\omega}\Upsilon(\bar{t},\tilde{\theta},\tilde{\Gamma}) \,.\label{eq:dotTildeGamma_20250603_1}
\end{align}
At this point, we can implement a suitable averaging mechanism within the transformed time scale $\bar{t}$ based on the dynamics (\ref{eq:dhatGdt_20250206_3})--(\ref{eq:dotTildeGamma_20250603_1}). Despite the non-periodicity of the triggering events and discontinuity on the right-hand sides, the closed-loop system maintains its periodicity over time due to the periodic probing and demodulation signals. This unique characteristic allows for the application of the averaging results established by Plotnikov \cite{P:1979} to this particular setup---see Appendix~\ref{appendix_plotnikov}. 

\subsection{Average Closed-Loop System}

Defining the augmented state as follows
\begin{align}
x(\bar{t}):=\begin{bmatrix} \hat{G}^{\top}(\bar{t})\,, \tilde{\theta}^{\top}(\bar{t})\,, \tilde{\Gamma}^{\top}(\bar{t})\end{bmatrix}^{\top}\,, \label{eq:X}
\end{align}
the system (\ref{eq:dhatGdt_20250206_3})--(\ref{eq:dotTildeGamma_20250603_1}) reduces to
\begin{align}
\dfrac{dx}{d\bar{t}}(\bar{t})&=\dfrac{1}{\omega}f\left(\bar{t},x\right)\,. \label{eq:dotX_event}
\end{align}
Note that (\ref{eq:dotX_event}) is characterized by a small parameter $1/\omega$ as well as a $T$-periodic function $f\left(\bar{t},x\right)$ in $\bar{t}$ and, thereby, the averaging method can be performed on  $f\left(\bar{t},x\right)$ at $\displaystyle \lim_{\omega\to \infty}\dfrac{1}{\omega}=0$, as shown in reference \cite{P:1979}. The averaging method allows for determining in what sense the behavior of a constructed average autonomous system approximates the behavior of the non-autonomous system (\ref{eq:dotX_event}). Of course, it can be inferred intuitively that in instances where the response of a system is significantly slower than its excitation, the response will predominantly be dictated by the average characteristics of the excitation. By employing the averaging technique to (\ref{eq:dotX_event}), we derive the following average system
\begin{align}
\dfrac{dx_{\rm{av}}}{d\bar{t}}(\bar{t})&=\dfrac{1}{\omega}f_{\rm{av}}\left(x_{\rm{av}}\right) \,, \label{eq:dotXav_event_1} \\
f_{\rm{av}}\left(x_{\rm{av}}\right)&=\dfrac{1}{T}\int_{0}^{T}f\left(\xi,x_{\rm{av}}\right)d\xi
\,.  \label{eq:mathcalFav_event}
\end{align}
Therefore, ``freezing'' the average states of $\hat{G}(\bar{t})$ and $\tilde{\theta}(\bar{t})$ in (\ref{eq:dhatGdt_20250206_3})--(\ref{eq:dtildeThetadt_20250206_3}), one gets for all $\bar t\in [\bar t_k, \bar t_{k+1})$, the nonlinear average system
\begin{align}
\frac{d\hat{G}_{\rm{av}}}{d\bar{t}}(\bar{t})&\mathbb{=} \mathbb{-}\frac{1}{\omega}H^{\ast}KH^{\ast-1}\hat{G}_{\rm{av}}(\bar{t})\mathbb{-}\frac{1}{\omega}H^{\ast}K\tilde{\Gamma}(\bar{t})\hat{G}_{\rm{av}}(\bar{t}) \nonumber \\
&\quad\mathbb{-}\frac{1}{\omega}H^{\ast}Ke_{\rm{av}}(\bar{t}), \label{eq:dhatGAv_20250603} \\
\frac{d\tilde{\theta}_{\rm{av}}}{d\bar{t}}(\bar{t})&\mathbb{=}\mathbb{-}\frac{1}{\omega}K\tilde{\theta}_{\rm{av}}(\bar{t})\mathbb{-}\frac{1}{\omega}K\tilde{\Gamma}_{\rm{av}}(\bar{t})H^{\ast}\tilde{\theta}_{\rm{av}}(\bar{t}) \mathbb{-}\frac{1}{\omega}Ke_{\rm{av}}(\bar{t}), \label{eq:dtildeThetaAv_20250603} \\
\frac{d\tilde{\Gamma}_{\rm{av}}}{d\bar{t}}(\bar{t})&=\!\mathbb{-}\frac{\omega_{r}}{\omega}\tilde{\Gamma}_{\rm{av}}(\bar{t})\mathbb{-}\frac{\omega_{r}}{\omega}\tilde{\Gamma}_{\rm{av}}(\bar{t})H^{\ast}\tilde{\Gamma}_{\rm{av}}(\bar{t})\,,\label{eq:dotTildeGammaAv_20250603} 
\end{align}
with 
\begin{align}
\hat{G}_{\rm{av}}(\bar{t})&= H^{\ast}\tilde{\theta}_{\rm{av}}(\bar{t})\,, \label{eq:hatGAv_20250603}  \\
\tilde{\Gamma}_{\rm{av}}(\bar{t})&=\Gamma_{\rm{av}}(\bar{t})-H^{\ast-1}\,. \label{eq:tildeGammaAv_20250603}  
\end{align}
Hence, by using (\ref{eq:hatGAv_20250603}) and (\ref{eq:tildeGammaAv_20250603}), 
\begin{align}
\Gamma_{\rm{av}}(\bar{t})\hat{G}_{\rm{av}}(\bar{t})&= \tilde{\theta}_{\rm{av}}(\bar{t})+\tilde{\Gamma}_{\rm{av}}(\bar{t})H^{\ast}\tilde{\theta}_{\rm{av}}(\bar{t})\,. \label{eq:GammaGAv_20250603}  
\end{align}

\subsection{Linearized Average Closed-Loop System}

Before analyze the stability properties of our static event-triggered multivariable Newton-based extremum seeking strategy, it is useful to find an average system representation of (\ref{eq:dhatGAv_20250603})--(\ref{eq:dotTildeGammaAv_20250603}) near the origin  $(\hat{G}_{\rm{av}}\,,\tilde{\theta}_{\rm{av}}\,,\tilde{\Gamma}_{\rm{av}})=(0\,,0\,,0)$. The averaged closed-loop dynamics of the system (\ref{eq:dhatGAv_20250603})--(\ref{eq:dotTildeGammaAv_20250603}), as well as (\ref{eq:GammaGAv_20250603}), have nonlinear terms that are quadratic in the $(\tilde{\Gamma}_{\rm{av}}\,,\hat{G}_{\rm{av}})$, $(\tilde{\Gamma}_{\rm{av}}\,,\tilde{\theta}_{\rm{av}})$ and $\tilde{\Gamma}_{\rm{av}}$, respectively: $H^{\ast}K\tilde{\Gamma}_{\rm{av}}(\bar{t})\hat{G}_{\rm{av}}(\bar{t})$, $\tilde{\Gamma}_{\rm{av}}(\bar{t})H^{\ast}\tilde{\theta}_{\rm{av}}(\bar{t})$ and $\tilde{\Gamma}_{\rm{av}}(\bar{t})H^{\ast}\tilde{\Gamma}_{\rm{av}}(\bar{t})$. These nonlinearities complicate the analysis of the averaged closed-loop system, but their influence diminishes near the origin, where the error and gradient magnitudes become small. 

Therefore, to simplify the analysis, the averaged closed-loop system  (\ref{eq:dhatGAv_20250603})--(\ref{eq:dotTildeGammaAv_20250603}) can be linearized around the equilibrium by 
\begin{align}
\frac{d\hat{G}_{\rm{av}}}{d\bar{t}}(\bar{t})&= -\frac{1}{\omega}H^{\ast}KH^{\ast-1}\hat{G}_{\rm{av}}(\bar{t})-\frac{1}{\omega}H^{\ast}Ke_{\rm{av}}(\bar{t}), \label{eq:dhatGAv_20250603_2} \\
\frac{d\tilde{\theta}_{\rm{av}}}{d\bar{t}}(\bar{t})&=-\frac{1}{\omega}K\tilde{\theta}_{\rm{av}}(\bar{t})-\frac{1}{\omega}Ke_{\rm{av}}(\bar{t}), \label{eq:dtildeThetaAv_20250603_2} \\
\frac{d\tilde{\Gamma}_{\rm{av}}}{d\bar{t}}(\bar{t})&=-\frac{\omega_{r}}{\omega}\tilde{\Gamma}_{\rm{av}}(\bar{t})\,,\label{eq:dotTildeGammaAv_20250603_2} 
\end{align}
and
\begin{align}
\Gamma_{\rm{av}}(\bar{t})\hat{G}_{\rm{av}}(\bar{t})&= \tilde{\theta}_{\rm{av}}(\bar{t})\,, \label{eq:GammaGAv_20250603_2}  
\end{align}
representing a local approximation of (\ref{eq:GammaGAv_20250603}).

From (\ref{eq:dhatGAv_20250603_2}) and (\ref{eq:dtildeThetaAv_20250603_2}), it is evident the Input-to-State Stability relationship of $\hat{G}_{\rm{av}}(\bar{t})$ and $\tilde{\theta}_{\rm{av}}(\bar{t})$ with respect to the average measurement error $e_{\rm{av}}(\bar{t})$. Moreover, from (\ref{eq:GammaGAv_20250603_2}), the average version of (\ref{eq:e_event_20250602_v1}) becomes simply
\begin{align}
    e_{\rm{av}}(\bar{t})&=\tilde{\theta}_{\rm{av}}(\bar{t}_{\kappa})-\tilde{\theta}_{\rm{av}}(\bar{t})\,. \label{eq:eAv_20250606}
\end{align}

Now, by using (\ref{eq:GammaGAv_20250603_2}), we introduce {\bf Definition~\ref{def:averageStaticEvent}} as an average version of {\bf Definition~\ref{def:staticEvent}}.

\begin{definition}[\small{Average Event-Triggering Condition}] \label{def:averageStaticEvent} The  average event-triggered condition consists of two components:
\begin{enumerate}
	\item A sequence of increasing times $\mathcal{I}$ such that $\mathcal{I}=\{\bar{t}_{0}\,, \bar{t}_{1}\,, \bar{t}_{2}\,,\ldots\}$ with $\bar{t}_{0}=0$, generated under the following rules:
		\begin{itemize}[leftmargin=-0.1in]
			\item If $\left\{\bar{t} \in\mathbb{R}^{+}: \bar{t} \mathbb{>} \bar{t}_{\kappa} ~ \mathbb{\wedge} ~ \sigma\|\tilde{\theta}_{\rm{av}}(\bar{t})\| \mathbb{-}\alpha\|e_{\rm{av}}(\bar{t})\| \mathbb{<} 0 \right\} \mathbb{=} \emptyset$, then the set of the times of the events is $\mathcal{I}=\{\bar{t}_{0}\,, \bar{t}_{1}\,, \ldots, \bar{t}_{\kappa}\}$.
			\item If $\left\{\bar{t} \in\mathbb{R}^{+}: \bar{t} \mathbb{>} \bar{t}_{\kappa} ~ \mathbb{\wedge} ~ \sigma\|\tilde{\theta}_{\rm{av}}(\bar{t})\| \mathbb{-}\alpha\|e_{\rm{av}}(\bar{t})\| \mathbb{<} 0 \right\} \mathbb{\neq} \emptyset$, the next event time, consisting of the static event-triggering mechanism, is given by
            \vspace{-0.3cm}
				\begin{align}
					\!\!\!\!\!\!\!\!\!\!\!\bar{t}_{\kappa+1}&\mathbb{=}\inf\left\{\bar{t} \in\mathbb{R}^{+}\mathbb{:} ~ \bar{t} \mathbb{>} \bar{t}_{\kappa} ~ \mathbb{\wedge} ~ \sigma\|\tilde{\theta}_{\rm{av}}(\bar{t})\| \mathbb{-}\alpha\|e_{\rm{av}}(\bar{t})\| \mathbb{<} 0 \right\}. \label{eq:tk+1_event_av}
				\end{align}
		\end{itemize}
        \vspace{-0.30cm}
	\item The feedback tuning action using the average estimation error $\tilde{\theta}_{\rm{av}}(\bar{t})$ updated at the triggering instants is given by
		\begin{align}
			u_{\kappa}^{\rm{av}}=-K\tilde{\theta}_{\rm{av}}(\bar{t}_{\kappa}) \,,  \quad \forall \bar{t} \in \lbrack \bar{t}_{\kappa}\,, \bar{t}_{\kappa+1}\phantom{(}\!\!)\,, \quad \kappa\in \mathbb{N}\,. \label{eq:U_MD2}
		\end{align}
\end{enumerate}  
\end{definition}

We claim that the event-triggering mechanism discussed above guarantee the asymptotic stabilization of $\hat{G}_{\rm{av}}(\bar{t})$ and, consequently, that of $\tilde{\theta}_{\rm{av}}(\bar{t})$. Moreover, since the matrix $H^{\ast}KH^{\ast-1}$ has the same eigenvalues as $K$, the convergence rate of both $\hat{G}_{\rm{av}}(\bar{t})$ and $\tilde{\theta}_{\rm{av}}(\bar{t})$ can be freely chosen by the designer.

Next, the stability analysis will be carried out considering the static event-triggering  mechanism introduced above. Note that, it is considered the total lack of knowledge of the nonlinear maps (\ref{eq:y_v2}), {\it i.e.}, $H^{\ast}$, $\theta^{\ast}$ and $Q^{\ast}$ are composed by completely unknown parameters.

\section{Stability Analysis}\label{ETESNC_unknownH*}

Theorem~\ref{thm:NETESC_2} states the local exponential practical stability of the extremum seeking control of Fig.~\ref{fig:BD_NewtonES_v5} based on event-triggering execution mechanism is ensured.

\begin{theorem} \label{thm:NETESC_2}
Consider the closed-loop average dynamics of the gradient estimate (\ref{eq:dhatGAv_20250603_2}), the average error vector (\ref{eq:eAv_20250606}), under  Assumptions \ref{assumption_w} and \ref{assumption_lyapunovEq}, and the average \textbf{static} event-triggering mechanism given by \textbf{Definition \ref{def:averageStaticEvent}}. For $\omega>0$ defined in (\ref{eq:omega_event_1}), sufficiently large, the equilibrium $(\hat{G}_{\rm{av}}\,,\tilde{\theta}_{\rm{av}}\,,\tilde{\Gamma}_{\rm{av}})=(0\,,0\,,0)$ is locally exponentially stable such that the following inequalities can be obtained for the non-average signals:
\begin{small}
\begin{align}
&\|\theta(t)-\theta^{\ast}\|\leq \sqrt{\frac{\lambda_{\max}(P_{1})}{\lambda_{\min}(P_{1})}} \times \nonumber \\
&\times \exp\left(\!\!\mathbb{-}\frac{1}{2}\frac{\lambda_{\min}(Q)}{\lambda_{\max}(P_{1})}(1\mathbb{-}\sigma)t\right)\|\theta(0)-\theta^{\ast}\|+\mathcal{O}\left(a+\frac{1}{\omega}\right) ,  \label{eq:normTheta_thm1} \\
    &|y(t) - Q^{\ast}| \leq 2\frac{\lambda_{\max}((-1)^{\rm{op}}H^{\ast})}{\lambda_{\min}((-1)^{\rm{op}}H^{\ast})}\frac{\lambda_{\max}(P_{1})}{\lambda_{\min}(P_{1})} \times \nonumber \\
    &\times \exp\left(-\frac{\lambda_{\min}(Q)}{\lambda_{\max}(P_{1})}\left(1-\sigma\right)t\right)|y(0)-Q^{\ast}| +\mathcal{O}\left(a^{2}+\frac{1}{\omega^{2}}\right), \label{eq:normY_thm1} \\
    &\|\hat{G}(t)\| \leq\sqrt{\frac{\lambda_{\max}(P_{1})}{\lambda_{\min}(P_{1})}}\|H^{\ast}\|\|H^{\ast-1}\| \times \nonumber \\
&\times\exp\left(-\frac{1}{2}\frac{\lambda_{\min}(Q)}{\lambda_{\max}(P_{1})}\left(1-\sigma\right)t\right)\|\hat{G}(0)\|+\mathcal{O}\left(\frac{1}{\omega}\right), \label{eq:normG_thm1} \\
    & \|\Gamma(t)-H^{\ast-1}\| \leq \|\exp\left(-\omega_{r} t\right)\|\|\Gamma(0)-H^{\ast-1}\|+\mathcal{O}\left(\frac{1}{\omega}\right). \label{eq:normGamma_thm1}
\end{align}
\end{small}
where $\rm{op} = 0$ if $H^{\ast} > 0$ and $\rm{op} = 1$ if $H^{\ast} < 0$, so that $(-1)^{\rm{op}} H^{\ast} > 0$ in both cases, and where $a=\sqrt{\sum_{i=1}^{n}a_{i}^{2}}$, with $a_i$ defined in (\ref{eq:S_v1}), (\ref{eq:M_v1}), (\ref{eq:Nii}) and (\ref{eq:Nij}). In addition, there exists a lower bound  $\tau^{\ast}$ for the inter-execution interval $t_{k+1}-t_{k}$  for all $k \in \mathbb{N}$ precluding the Zeno behavior.
\end{theorem} 
\textit{Proof:} The proof of the theorem is divided into two parts: stability analysis and avoidance of Zeno behavior.

\begin{flushleft}
\underline{\it A. Stability Analysis}
\end{flushleft}

Now, consider the following Lyapunov function candidate for the average system (\ref{eq:dtildeThetaAv_20250603_2}):
\begin{align}
V_{\rm{av}}(\bar{t})=\tilde{\theta}^{\top}_{\rm{av}}(\bar{t})P_{1}\tilde{\theta}_{\rm{av}}(\bar{t}) \,,\, \quad P_{1}=P_{1}^{\top} >0, \label{eq:Vav_event_pf2}
\end{align}  
with time-derivative
\begin{align}
\frac{dV_{\rm{av}}(\bar{t})}{d\bar{t}}&=-\frac{1}{\omega}\tilde{\theta}_{\rm{av}}^{\top}(\bar{t})Q\tilde{\theta}_{\rm{av}}(\bar{t})+\frac{2}{\omega}\tilde{\theta}_{\rm{av}}^{\top}(\bar{t})P_{1}KH^{\ast}e_{\rm{av}}(\bar{t})\,, \label{eq:dotVav_event_1_pf2}
\end{align}
whose  upper bound satisfies
\begin{align}
\frac{dV_{\rm{av}}(\bar{t})}{d\bar{t}}&\leq\mathbb{-}\frac{\lambda_{\min}(Q)}{\omega}\|\tilde{\theta}_{\rm{av}}(\bar{t})\|^{2} \nonumber \\
&\quad~\mathbb{+}
\frac{ 2\|P_{1}KH^{\ast}\|}{\omega} \|\tilde{\theta}_{\rm{av}}(\bar{t})\|\|e_{\rm{av}}(\bar{t})\|. \label{eq:dotVav_event_2_pf2}
\end{align}

In the proposed event-triggering mechanism, the update law is (\ref{eq:tk+1_event_av}), the vector $u_{\rm{av}}(\bar{t})$ is held constant between two consecutive events, and therefore the norm of average measurement error $e_{\rm{av}}(\bar{t})$ is upper bounded by
\begin{align}
\|e_{\rm{av}}(\bar{t})\|& \leq \frac{\sigma}{\alpha}\|\tilde{\theta}_{\rm{av}}(\bar{t})\|\,, \quad \sigma \in (0,1)\,, \quad \alpha>\frac{ 2\|P_{1}KH^{\ast}\|}{\lambda_{\min}(Q)}  \,. \label{eq:eAv_upperBound}
\end{align}
Now, by using (\ref{eq:eAv_upperBound}), inequality (\ref{eq:dotVav_event_2_pf2}) is upper bounded as 
\begin{align}
\frac{dV_{\rm{av}}(\bar{t})}{d\bar{t}}&\leq-\frac{\lambda_{\min}(Q)}{\omega}\left(1-\sigma\right)\|\tilde{\theta}_{\rm{av}}(\bar{t})\|^{2}\,. \label{eq:dotVav_event_4_pf2}
\end{align}

By using the Rayleigh-Ritz Inequality \cite{K:2002}, 
\begin{align}
\lambda_{\min}(P_{1})\|\tilde{\theta}_{\rm{av}}(\bar{t})\|^{2}\leq V_{\rm{av}}(\bar{t}) \leq \lambda_{\max}(P_{1})\|\tilde{\theta}_{\rm{av}}(\bar{t})\|^{2}\,, \label{eq:Rayleigh-Ritz_pf2}
\end{align}
 and the following upper bound for (\ref{eq:dotVav_event_4_pf2}), one obtains
\begin{align}
\frac{dV_{\rm{av}}(\bar{t})}{d\bar{t}}&\leq-\frac{1}{\omega}\frac{\lambda_{\min}(Q)}{\lambda_{\max}(P_{1})}\left(1-\sigma\right)V_{\rm{av}}(\bar{t}) \,. \label{eq:dotVav_event_5_pf2}
\end{align}

Then, by invoking the Comparison Principle \cite[Comparison Lemma]{K:2002}, an upper bound $\bar{V}_{\rm{av}}(\bar{t})$ for $V_{\rm{av}}(\bar{t})$ can be derived
\begin{align}
V_{\rm{av}}(\bar{t})\leq \bar{V}_{\rm{av}}(\bar{t}) \,, \quad \forall \bar{t}\in \lbrack \bar{t}_{k},\bar{t}_{k+1}\phantom{(}\!\!) \,, \label{eq:VavBarVav_1_pf2}
\end{align}
 is given by the solution of:
\begin{align}
\frac{d\bar{V}_{\rm{av}}(\bar{t})}{d\bar{t}}\mathbb{=}\mathbb{-}\frac{1}{\omega}\frac{\lambda_{\min}(Q)}{\lambda_{\max}(P_{1})}\left(1\mathbb{-}\sigma\right)\bar{V}_{\rm{av}}(\bar{t})\,,~ \bar{V}_{\rm{av}}(\bar{t}_{k})\mathbb{=}V_{\rm{av}}(\bar{t}_{k})\,.
\end{align}
In other words, $ \forall \bar{t}\in \lbrack \bar{t}_{k},\bar{t}_{k+1}\phantom{(}\!\!)$,
\begin{align}
\bar{V}_{\rm{av}}(\bar{t})=\exp\left(\mathbb{-}\frac{1}{\omega}\frac{\lambda_{\min}(Q)}{\lambda_{\max}(P_{1})}\left(1\mathbb{-}\sigma\right)\bar{t}\right)V_{\rm{av}}(\bar{t}_{k})\,, \label{eq:_pf2}
\end{align}
and inequality (\ref{eq:VavBarVav_1_pf2}) is rewritten, $\forall \bar{t}\in \lbrack \bar{t}_{k},\bar{t}_{k+1}\phantom{(}\!\!)$, as
\begin{align}
V_{\rm{av}}(\bar{t})\leq \exp\left(\mathbb{-}\frac{1}{\omega}\frac{\lambda_{\min}(Q)}{\lambda_{\max}(P_{1})}\left(1\mathbb{-}\sigma\right)\bar{t}\right)V_{\rm{av}}(\bar{t}_{k})
\,. \label{eq:VavBarVav_2_pf2}
\end{align}

By defining $\bar{t}_{k}^{+}$ and $\bar{t}_{k}^{-}$ as the right and left limits of $\bar{t}=\bar{t}_{k}$, respectively, it easy to verify that $V_{\rm{av}}(\bar{t}_{k+1}^{-})\leq \exp\left(\mathbb{-}\frac{1}{\omega}\frac{\lambda_{\min}(Q)}{\lambda_{\max}(P_{1})}\left(1\mathbb{-}\sigma\right)(\bar{t}_{k+1}^{-}-\bar{t}_{k}^{+})\right)V_{\rm{av}}(\bar{t}_{k}^{+})$. Since $V_{\rm{av}}(\bar{t})$ is continuous, $V_{\rm{av}}(\bar{t}_{k+1}^{-})=V_{\rm{av}}(\bar{t}_{k+1})$ and $V_{\rm{av}}(\bar{t}_{k}^{+})=V_{\rm{av}}(\bar{t}_{k})$, thus
\begin{align}
    V_{\rm{av}}(\bar{t}_{k+1})\mathbb{\leq} \exp\left(\!\!\mathbb{-}\frac{1}{\omega}\frac{\lambda_{\min}(Q)}{\lambda_{\max}(P_{1})}\left(1\mathbb{-}\sigma\right)(\bar{t}_{k\mathbb{+}1}\mathbb{-}\bar{t}_{k})\!\!\right)\! V_{\rm{av}}(\bar{t}_{k}). \label{METES_eq:mmd_1_s}
\end{align}

The recursive structure of (\ref{METES_eq:mmd_1_s}) allows us to derive an upper bound for the evolution of the Lyapunov function (\ref{eq:Vav_event_pf2}) over successive intervals such that we can iteratively relate the Lyapunov function at the current time $\bar{t}$ to its value at earlier triggering times. This approach leads to a product of exponential terms, each corresponding to an inter-execution interval, with decay rates determined by the system parameters \(\omega\), \(\sigma\), and the spectral properties of \(P_1\) and \(Q\). By carefully reorganizing these terms and using the fact that the sum of all time intervals up to \(\bar{t}\) equals \(\bar{t}\) itself, we obtain a compact expression for global exponential convergence from the initial condition, {\it i.e.}, for any $\bar{t}\geq 0$ in $ \bar{t}\in \lbrack \bar{t}_{k},\bar{t}_{k+1}\phantom{(}\!\!)$, $k \in \mathbb{N}$, one has 
\begin{align}
    &V_{\rm{av}}(\bar{t})\leq \exp\left(\!\!\mathbb{-}\frac{1}{\omega}\frac{\lambda_{\min}(Q)}{\lambda_{\max}(P_{1})}\left(1\mathbb{-}\sigma\right)(\bar{t}-\bar{t}_{k})\right) V_{\rm{av}}(\bar{t}_{k}) \nonumber \\
    &\leq \exp\left(\!\!\mathbb{-}\frac{1}{\omega}\frac{\lambda_{\min}(Q)}{\lambda_{\max}(P_{1})}\left(1\mathbb{-}\sigma\right)(\bar{t}-\bar{t}_{k})\right)  \nonumber \\
		&\quad \times\exp\left(\!\!\mathbb{-}\frac{1}{\omega}\frac{\lambda_{\min}(Q)}{\lambda_{\max}(P_{1})}\left(1\mathbb{-}\sigma\right)(\bar{t}_{k}-\bar{t}_{k-1})\right)V_{\rm{av}}(\bar{t}_{k-1}) \nonumber \\
    &\leq \ldots \leq \nonumber \\
    &\leq \exp\left(\!\!\mathbb{-}\frac{1}{\omega}\frac{\lambda_{\min}(Q)}{\lambda_{\max}(P_{1})}\left(1\mathbb{-}\sigma\right)(\bar{t}\mathbb{-}\bar{t}_{k})\right) \nonumber \\
		&\quad \times \prod_{i=1}^{i=k}\exp\left(\!\!\mathbb{-}\frac{1}{\omega}\frac{\lambda_{\min}(Q)}{\lambda_{\max}(P_{1})}\left(1\mathbb{-}\sigma\right)(\bar{t}_{i}-\bar{t}_{i-1})\right)V_{\rm{av}}(\bar{t}_{i-1}) \nonumber \\
    &=\exp\left(\!\!\mathbb{-}\frac{1}{\omega}\frac{\lambda_{\min}(Q)}{\lambda_{\max}(P_{1})}\left(1\mathbb{-}\sigma\right)\bar{t}\right) V_{\rm{av}}(0)\,, \quad \forall \bar{t}\geq 0\,. \label{METES_eq:VavBarVav_2_pf2}
\end{align}

Now, lower bounding the left-hand side and upper bounding the right-hand size of (\ref{METES_eq:VavBarVav_2_pf2}) with the corresponding sides of (\ref{eq:Rayleigh-Ritz_pf2}), one gets
\begin{align}
&\lambda_{\mathbb{\min}}(P_{1})\|\tilde{\theta}_{\rm{av}}(\bar{t})\|^{2}\leq\nonumber \\
&\exp\left(-\frac{1}{\omega}\frac{\lambda_{\min}(Q)}{\lambda_{\max}(P_{1})}\left(1-\sigma\right)\bar{t}\right)\lambda_{\max}(P_{1})\|\tilde{\theta}_{\rm{av}}(0)\|^{2}\,. \label{eq:VavBarVav_3_pf2}
\end{align}
Therefore,
\begin{align}
&\|\tilde{\theta}_{\rm{av}}(\bar{t})\|\leq \nonumber \\
&\exp\left(-\frac{1}{2\omega}\frac{\lambda_{\min}(Q)}{\lambda_{\max}(P_{1})}\left(1-\sigma\right)\bar{t}\right)\sqrt{\frac{\lambda_{\max}(P_{1})}{\lambda_{\min}(P_{1})}}\|\tilde{\theta}_{\rm{av}}(0)\|\,. \label{eq:normHatTildeThetaAv_20250616}
\end{align}
Since (\ref{eq:dtildeThetadt_20250206_3}) has a discontinuous right-hand side, but it is also $T$-periodic in $t$, and noting that the average system with state $\tilde{\theta}_{\rm{av}}(\bar{t})$ is asymptotically stable according to (\ref{eq:normHatTildeThetaAv_20250616}), we can invoke the averaging theorem in \cite[Theorem~2]{P:1979} (see also Appendix~\ref{appendix_plotnikov}) to conclude that
\begin{align}
\|\tilde{\theta}(t)-\tilde{\theta}_{\rm{av}}(t)\|\leq\mathcal{O}\left(\frac{1}{\omega}\right)\,. \label{eq:plotnikov}
\end{align}
Now, adding and subtracting $\tilde{\theta}_{\rm{av}}(\bar{t})$ in the right-hand side of vector form of (\ref{eq:thetai_v1}), in the time scale $\bar{t}$, one has
\begin{align}
\theta(\bar{t}) - \theta^\ast = \tilde{\theta}_{\rm{av}}(\bar{t})+\tilde{\theta}(\bar{t})-\tilde{\theta}_{\rm{av}}(\bar{t}) + S(\bar{t})\,, \label{eq:q_ETSSC_v3cu}
\end{align}  
whose norm can be upper bounded by using the triangle inequality \cite{A:1957}, one gets
\begin{align}
\|\theta(\bar{t}) - \theta^\ast\| \leq \|\tilde{\theta}_{\rm{av}}(\bar{t})\|+\|\tilde{\theta}(\bar{t})-\tilde{\theta}_{\rm{av}}(\bar{t})\| + \|S(\bar{t})\|\,.\label{eq:q_ETSSC_v3cu2}
\end{align} 
Since sine functions in (\ref{eq:S_v1}) are uniformly bounded, it is possible to derive a uniform upper bound for the Euclidean norm of $S(\bar{t})$ such that
\begin{align}
    \|S(\bar{t})\|\leq \sqrt{\sum_{i=1}^{n}a_{i}^{2}}=a \,.\label{eq:S_v2}
\end{align}
Thus, by using (\ref{eq:normHatTildeThetaAv_20250616}), (\ref{eq:plotnikov}) and (\ref{eq:S_v2}), inequality (\ref{eq:q_ETSSC_v3cu2}) is upper bounded by 
\begin{align}
&\|\theta(\bar{t})-\theta^{\ast}\|  \leq  \nonumber \\
& \exp\left(-\frac{1}{2\omega}\frac{\lambda_{\min}(Q)}{\lambda_{\max}(P_{1})}\left(1-\sigma\right)\bar{t}\right)\sqrt{\frac{\lambda_{\max}(P_{1})}{\lambda_{\min}(P_{1})}}\|\theta(0)-\theta^{\ast}\| \nonumber \\
&+\mathcal{O}\left(a+\frac{1}{\omega}\right)\,.\label{eq:q_ETSSC_v4}
\end{align} 
Therefore, by using (\ref{eq:q_ETSSC_v4}) and the original system time scale $t=\dfrac{\bar{t}}{\omega}$, inequality (\ref{eq:normTheta_thm1}) is verified.

Notice that the Hessian matrix $H^{\ast} \in \mathbb{R}^{n \times n}$ in (\ref{eq:y_v2}) is symmetric and sign-definite ($H^{\ast} > 0$ or $H^{\ast} < 0$) such that the quadratic term $(\theta(t)-\theta^{\ast})^{\top}H^{\ast}(\theta(t)-\theta^{\ast})$ is bounded above and below by quadratic functions of $\|\theta(t)-\theta^{\ast}\|$ with proportional constants given by the minimum and maximum eigenvalues of $H^{\ast}$. From (\ref{eq:y_v2}), the term $|y(t) - Q^{\ast}|$ corresponds to the magnitude of a quadratic form, which can be expressed as $\dfrac{1}{2} |(\theta(t)-\theta^{\ast})^{\top}H^{\ast}(\theta(t)-\theta^{\ast})|$. To uniformly treat both cases where $H^{\ast}$ is positive or negative definite, we introduce the factor $\rm{op} =\begin{cases} 0 \,, \quad H^{\ast}>0 \\ 1 \,, \quad H^{\ast}<0 \end{cases}$, so that $(-1)^{\rm{op}} H^{\ast} > 0$ in both cases. Consequently, $|y(t) - Q^{\ast}| = \dfrac{1}{2} |(\theta(t)-\theta^{\ast})^{\top}H^{\ast}(\theta(t)-\theta^{\ast})|$ can be bounded by $ \frac{\lambda_{\min}((-1)^{\rm{op}}H^{\ast})}{2}\|\theta(t)-\theta^{\ast}\|^2 \leq |y(t) - Q^{\ast}| \leq \frac{\lambda_{\max}((-1)^{\rm{op}}H^{\ast})}{2}\|\theta(t)-\theta^{\ast}\|^2$, being valid uniformly, regardless of the sign of $H^{\ast}$. Thus, by using (\ref{eq:q_ETSSC_v4}) and Young's inequality \cite{K:2002}, the following upper bound can verified 
\begin{align}
    &|y(\bar{t}) - Q^{\ast}| \leq 2\frac{\lambda_{\max}((-1)^{\rm{op}}H^{\ast})}{\lambda_{\min}((-1)^{\rm{op}}H^{\ast})}\frac{\lambda_{\max}(P_{1})}{\lambda_{\min}(P_{1})}  \nonumber \\
    &\times \exp\left(-\frac{1}{\omega}\frac{\lambda_{\min}(Q)}{\lambda_{\max}(P_{1})}\left(1-\sigma\right)\bar{t}\right)|y(0)-Q^{\ast}| \nonumber \\
&+\mathcal{O}\left(a^{2}+\frac{1}{\omega^{2}}\right)\,,
\end{align}
such that, in the time scale $t$, inequality (\ref{eq:normY_thm1}) holds.   

Since $H^{\ast}$ is invertible, from (\ref{eq:hatGAv_20250603}), it is possible to verify that $\|H^{\ast}\|^{-1} \|\hat{G}_{\rm{av}}(\bar{t})\| \leq \|\tilde{\theta}_{\rm{av}}(\bar{t})\| \leq \|H^{\ast-1}\| \|\hat{G}_{\rm{av}}(\bar{t})\|$. Consequently, from (\ref{eq:normHatTildeThetaAv_20250616}), one has
\begin{align}
\|\hat{G}_{\rm{av}}(\bar{t})\|&\leq \exp\left(-\frac{1}{2\omega}\frac{\lambda_{\min}(Q)}{\lambda_{\max}(P_{1})}\left(1-\sigma\right)\bar{t}\right) \times \nonumber \\
&\quad\times\sqrt{\frac{\lambda_{\max}(P_{1})}{\lambda_{\min}(P_{1})}}\|H^{\ast}\|\|H^{\ast-1}\|\|\hat{G}_{\rm{av}}(0)\|\,. \label{eq:normHatGav_1_pf2}
\end{align}
Since the right-hand side of (\ref{eq:dhatGdt_20250206_3}) is both discontinuous and periodic in time with period $T$, and considering that the corresponding averaged system—governed by the state variable $\hat{G}_{\rm av}(\bar{t})$--is exponentially stable as established in (\ref{eq:normHatGav_1_pf2}), we can apply the averaging result from~\cite[Theorem~2]{P:1979} to conclude that
\begin{align}
\|\hat{G}(t)-\hat{G}_{\rm{av}}(t)\|\leq\mathcal{O}\left(\frac{1}{\omega}\right)\,. \label{eq:plotnikovG}
\end{align}
Consequently, by using (\ref{eq:normHatGav_1_pf2}) and (\ref{eq:plotnikovG}), we also can get
\begin{align}
&\|\hat{G}(\bar{t})\|=\|\hat{G}_{\rm{av}}(\bar{t})+\hat{G}(\bar{t})-\hat{G}_{\rm{av}}(\bar{t})\| \nonumber \\
&\leq\|\hat{G}_{\rm{av}}(\bar{t})\|+\|\hat{G}(\bar{t})-\hat{G}_{\rm{av}}(\bar{t})\| \nonumber \\
&\leq\sqrt{\frac{\lambda_{\max}(P_{1})}{\lambda_{\min}(P_{1})}}\exp\left(-\frac{1}{2\omega}\frac{\lambda_{\min}(Q)}{\lambda_{\max}(P_{1})}\left(1-\sigma\right)\bar{t}\right) \times \nonumber \\
&\quad\times\|H^{\ast}\|\|H^{\ast-1}\|\|\hat{G}(0)\|+\mathcal{O}\left(\frac{1}{\omega}\right)\,. \label{eq:normG_20250620}
\end{align}
Thus, in the time scale $t$, inequality (\ref{eq:normG_thm1}) holds.   

Analogously, recalling that the differential equation~\eqref{eq:dhatGdt_20250206_3} has a continuous right-hand side and is $T$-periodic in time, and meting that the solution of (\ref{eq:dotTildeGamma_20250603_1}) is given by 
\begin{align}
\tilde{\Gamma}_{\rm{av}}(\bar{t}) = \exp\left(-\frac{\omega_{r}}{\omega}\bar{t}\right)\tilde{\Gamma}_{\rm{av}}(0)\,, \label{eq:GammaTilde_20250623}
\end{align}
the conditions required to apply the classical averaging result stated in \cite[Chapter~10]{K:2002} are satisfied. Consequently, 
\begin{align}
\|\tilde{\Gamma}(\bar{t})-\tilde{\Gamma}_{\rm{av}}(\bar{t})\|\leq\mathcal{O}\left(\frac{1}{\omega}\right)\,, \label{eq:plotnikovGammaTilde}
\end{align}
and the behavior of the original time-varying system can be inferred from the stability properties of its averaged counterpart such that
\begin{align}
    \|\tilde{\Gamma}(\bar{t})\|&=\|\tilde{\Gamma}_{\rm{av}}(\bar{t})+\tilde{\Gamma}(\bar{t})-\tilde{\Gamma}_{\rm{av}}(\bar{t})\| \nonumber \\
    &\leq \|\tilde{\Gamma}_{\rm{av}}(\bar{t})\|+\|\tilde{\Gamma}(\bar{t})-\tilde{\Gamma}_{\rm{av}}(\bar{t})\| \nonumber \\
    &\leq\left \|\exp\left(-\frac{\omega_{r}}{\omega}\bar{t}\right)\right\|\|\tilde{\Gamma}_{\rm{av}}(0)\|+\mathcal{O}\left(\frac{1}{\omega}\right)\,.
\end{align}
Therefore, in the time scale $t$, inequality (\ref{eq:normGamma_thm1}) holds. 

\begin{flushleft}
\textcolor{black}{\underline{\it B. Avoidance of Zeno Behavior}}
\end{flushleft}

Since the average closed-loop system considering the sub-system of $\tilde{\theta}(\bar{t})$ consists of (\ref{eq:dtildeThetaAv_20250603_2}), with the event-triggering mechanism  (\ref{eq:tk+1_event_av}) and the average control law (\ref{eq:U_MD2}), we can conclude that $\|e_{\rm{av}}(\bar{t})\| \leq \bar{\sigma}\|\tilde{\theta}_{\rm{av}}(\bar{t})\|$ where $\bar{\sigma}=\dfrac{\sigma}{\alpha}$, resulting in, $ \forall \bar{t}\in \lbrack \bar{t}_{k},\bar{t}_{k+1}\phantom{(}\!\!)$,
\begin{align}
\bar{\sigma}\|\tilde{\theta}_{\rm{av}}(\bar{t})\|^{2}-\|e_{\rm{av}}(\bar{t})\|\|\tilde{\theta}_{\rm{av}}(\bar{t})\|\geq 0\,. \label{ineq:interEvents_1_static}
\end{align}
Using the Peter-Paul inequality \cite{W:1971}, $cd\leq \frac{c^2}{2\epsilon}+\frac{\epsilon d^2}{2}$ for all $c,d,\epsilon>0$, with $c=\|e_{\rm{av}}(\bar{t})\|$, $d=\|\tilde{\theta}_{\rm{av}}(\bar{t})\|$ and $\epsilon=\bar{\sigma}$, the inequality (\ref{ineq:interEvents_1_static}) is lower bounded by
\begin{align}
&\bar{\sigma} \|\tilde{\theta}_{\rm{av}}(\bar{t})\|^{2}-\|e_{\rm{av}}(\bar{t})\|\|\tilde{\theta}_{\rm{av}}(\bar{t})\|\geq \nonumber \\
&\bar{\sigma} \|\tilde{\theta}_{\rm{av}}(\bar{t})\|^{2}-\left(\frac{\bar{\sigma}}{2}\|\tilde{\theta}_{\rm{av}}(\bar{t})\|^2+\frac{1}{2\bar{\sigma}}\|e_{\rm{av}}(\bar{t})\|^2\right)\nonumber \\
&= q\|\tilde{\theta}_{\rm{av}}(\bar{t})\|^{2}-p\|e_{\rm{av}}(\bar{t})\|^2\,,\label{ineq:interEvents_2_static_pf2}
\end{align}
where 
\begin{align}
q&=\frac{\bar{\sigma}}{2} \quad \mbox{and} \quad p=\frac{1}{2\bar{\sigma}}\,.\label{ineq:interEvents_3_static_pf2}
\end{align} 
In \cite{G:2014}, it is shown that a lower bound for the inter-execution interval is given by the time duration it takes for the function
\begin{align}
\phi_{\rm{av}}(\bar{t})=\sqrt{\frac{p}{q}}\frac{\|e_{\rm{av}}(\bar{t})\|}{\|\tilde{\theta}_{\rm{av}}(\bar{t})\|} \label{eq:phi_1_static_pf2}
\end{align}
to go from 0 to 1. The time-derivative of (\ref{eq:phi_1_static_pf2}) is 
\begin{align}
\frac{d\phi_{\rm{av}}(\bar{t})}{d\bar{t}}&=\sqrt{\frac{p}{q}}\frac{1}{\|e_{\rm{av}}(\bar{t})\|\|\tilde{\theta}_{\rm{av}}(\bar{t})\|}\left[e_{\rm{av}}^{\top}(\bar{t})\frac{de_{\rm{av}}(\bar{t})}{d\bar{t}} \right. \nonumber \\
&\quad\left.-\tilde{\theta}_{\rm{av}}^{\top}(\bar{t})\frac{d\tilde{\theta}_{\rm{av}}(\bar{t})}{d\bar{t}}\left(\frac{\|e_{\rm{av}}(\bar{t})\|}{\|\tilde{\theta}_{\rm{av}}(\bar{t})\|}\right)^2\right]\,. \label{eq:dotPhi_1_static_pf2}
\end{align}
Now, by using equations (\ref{eq:dtildeThetaAv_20250603_2}) and (\ref{eq:eAv_20250606}) into (\ref{eq:dotPhi_1_static_pf2}), one arrives to 
\begin{align}
&\frac{d\phi_{\rm{av}}(\bar{t})}{d\bar{t}}\mathbb{=}\frac{1}{\omega}\sqrt{\frac{p}{q}}\frac{1}{\|e_{\rm{av}}(\bar{t})\|\|\tilde{\theta}_{\rm{av}}(\bar{t})\|} \nonumber \\
& \mathbb{\times}\Biggl\{\mathbb{-}e_{\rm{av}}^{\top}(\bar{t})Ke_{\rm{av}}(\bar{t}) \mathbb{-}e_{\rm{av}}^{\top}(\bar{t})K\tilde{\theta}_{\rm{av}}(\bar{t})+\Biggr.\nonumber \\
&\left. \mathbb{-}\left[\tilde{\theta}_{\rm{av}}^{\top}(\bar{t})K\tilde{\theta}_{\rm{av}}(\bar{t})\mathbb{+}\tilde{\theta}_{\rm{av}}^{\top}(\bar{t})Ke_{\rm{av}}(\bar{t})\right]\!\!\left(\frac{\|e_{\rm{av}}(\bar{t})\|}{\|\tilde{\theta}_{\rm{av}}(\bar{t})\|}\right)^{\!\!\! 2}\!\right\}\!\!. \label{eq:dotPhi_1_1_static_pf2}
\end{align}
Then, the following estimate holds:
\begin{align}
\frac{d\phi_{\rm{av}}(\bar{t})}{d\bar{t}}&\leq
\frac{\|K\|}{\omega}\sqrt{\frac{p}{q}}\left(1+\frac{\|e_{\rm{av}}(\bar{t})\|}{\|\tilde{\theta}_{\rm{av}}(\bar{t})\|}\right)^2\,. \label{eq:dotPhi_2_static_pf2}
\end{align}

Hence,  using (\ref{eq:phi_1_static_pf2}), inequality (\ref{eq:dotPhi_2_static_pf2}) is rewritten as
\begin{align}
\frac{d\phi_{\rm{av}}(\bar{t})}{d\bar{t}}&\leq\frac{\|K\|}{\omega}\sqrt{\frac{q}{p}}\left(\sqrt{\frac{p}{q}}+\phi_{\rm{av}}(\bar{t})\right)^{2}\,. \label{eq:dotPhi_3_static_pf2}
\end{align}
By invoking \cite[Comparison Lemma]{K:2002}, an upper bound $\tilde{\phi}_{\rm{av}}(\bar{t})$ for $\phi_{\rm{av}}(\bar{t})$ according to 
\begin{align}
\phi_{\rm{av}}(\bar{t})\leq \tilde{\phi}_{\rm{av}}(\bar{t}) \,, \quad \phi_{\rm{av}}(0)= \tilde{\phi}_{\rm{av}}(0)=0 \,, \label{eq:tildePhi_v2}
\end{align}
is given by the solution of the equation
\begin{align}
\frac{d\tilde{\phi}\textcolor{black}{_{\rm{av}}}(\bar{t})}{d\bar{t}}&=\frac{\|K\|}{\omega}\sqrt{\frac{q}{p}}\left(\sqrt{\frac{p}{q}}+\tilde{\phi}_{\rm{av}}(\bar{t})\right)^{2}\,. \label{eq:dotTildePhi_v2}
\end{align}
The solution of (\ref{eq:dotTildePhi_v2}), with the initial condition $\tilde{\phi}_{\rm{av}}(0) = 0$, is given by
\begin{align}
\tilde{\phi}_{\rm{av}}(\bar{t}) = \dfrac{\sqrt{\dfrac{p}{q}}}{1 - \dfrac{\|KH\|}{\omega} \dfrac{q}{p} \bar{t}} - \sqrt{\frac{p}{q}}.  \label{eq:phi-bart_v2}
\end{align}
Since $\phi_{\rm{av}}(t)$ in (\ref{eq:phi_1_static_pf2}) is an average version of $\phi(t)=\sqrt{\frac{p}{q}}\frac{\|e(\bar{t})\|}{\|\Gamma(\bar{t})\hat{G}(\bar{t})\|}$, by invoking \cite[Theorem 2]{P:1979}, one can write  
\begin{align}
|\phi(t)-\tilde{\phi}_{\rm{av}}(t)|\leq\mathcal{O}\left(\frac{1}{\omega}\right)\,.
\end{align}
By using the Triangle inequality \cite{A:1957}, one has
\begin{align}
\phi(t)&\leq\ \phi_{\rm{av}}(t)+\mathcal{O}\left(\frac{1}{\omega}\right) \leq \tilde{\phi}_{\rm{av}}(t)+\mathcal{O}\left(\frac{1}{\omega}\right)\! \nonumber \\
&=\dfrac{\sqrt{\frac{p}{q}}}{1 - \dfrac{\|K\|}{\omega} \dfrac{q}{p} t} - \sqrt{\dfrac{p}{q}} +\mathcal{O}\left(\frac{1}{\omega}\right). \label{eq:phi_t}
\end{align}
Now, defining 
\begin{align}
\hat{\phi}(t):=\dfrac{\sqrt{\dfrac{p}{q}}}{1 - \|K\| \dfrac{q}{p} t} - \sqrt{\dfrac{p}{q}} +\mathcal{O}\left(\frac{1}{\omega}\right)\,, \label{eq:hatPhi_t}
\end{align}
a lower bound on the inter-execution interval of the original system using the static event-triggered multivariable Newton-based extremum seeking is given by the time it takes for the function (\ref{eq:hatPhi_t}) to go from 0 to 1.  This is at least equal to 
\begin{align} \label{cataflan_novalgina}
\tau^{\ast}&=\dfrac{1}{\|K\|}\frac{1}{\bar{\sigma}^2}\frac{1-\mathcal{O}(1/\omega)}{1+1/\bar{\sigma}-\mathcal{O}(1/\omega)}\,,
\end{align} 
and the Zeno behavior is avoided in the original system as well. \hfill $\square$

\section{Simulation Results}

To evaluate and discuss the performance of continuous-time and static event-triggered implementations of the multivariable Newton-based extremum seeking control, we consider the nonlinear map defined in (\ref{eq:y_v2}), with a two-dimensional input \(\theta(t) \in \mathbb{R}^2\) and a scalar output \(y(t) \in \mathbb{R}\). The optimization problem is thus posed with the unknown parameters:
$$H = \begin{bmatrix} 100 & 30 \\ 30 & 20 \end{bmatrix} > 0,\quad Q^{\ast} = 100,\quad \theta^{\ast} = \begin{bmatrix} 2 \\ 4 \end{bmatrix}.$$
The perturbation and demodulation signals are configured with amplitudes \(a_1 = a_2 = 0.1\) and excitation frequencies \(\omega_1 = 1\,\text{rad/s}\), \(\omega_2 = 7\,\text{rad/s}\), following the design principles established in \cite{GKN:2012}. The static event-triggered mechanism is parameterized by \(\sigma = 0.75\) and \(\alpha = 0.8\). The Newton-based controller uses the gain matrix $$K = \begin{bmatrix} 1 & 0 \\ 0 & 1 \end{bmatrix}.$$ Both control strategies are initialized from the same initial parameter estimate: \(\hat{\theta}(0) = \begin{bmatrix} 2.5 \\ 5 \end{bmatrix}\).

\begin{figure}[h!]
	\centering
	\subfigure[\underline{\bf ET-Gradient-ES:} input of the map, $\theta(t)$. \label{fig:ETGradientES_theta}]{\includegraphics[width=4.1cm]{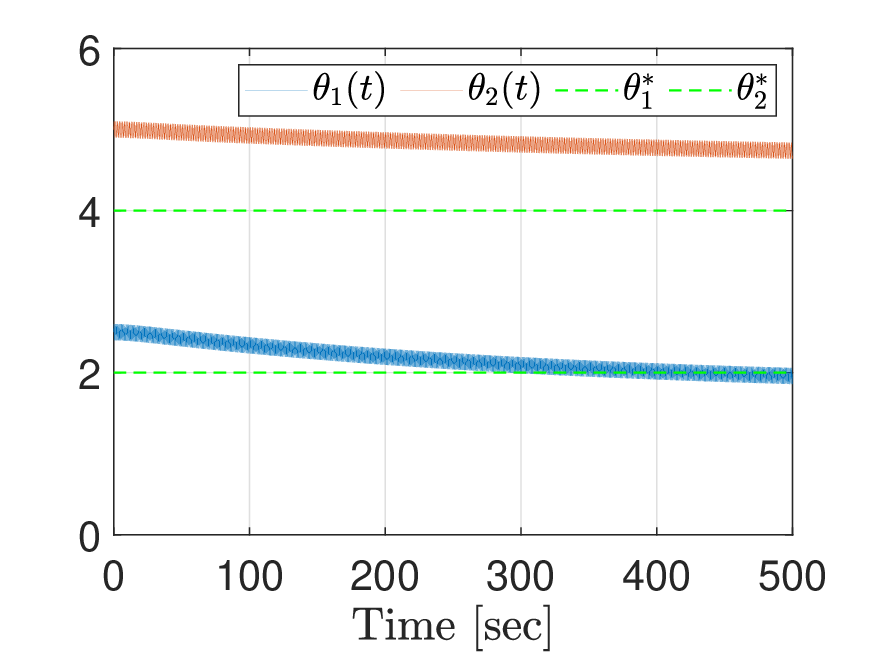}}
	~
	\subfigure[\underline{\bf ET-Newton-ES:} input of the map, $\theta(t)$. \label{fig:ETNewtonES_theta}]{\includegraphics[width=4.1cm]{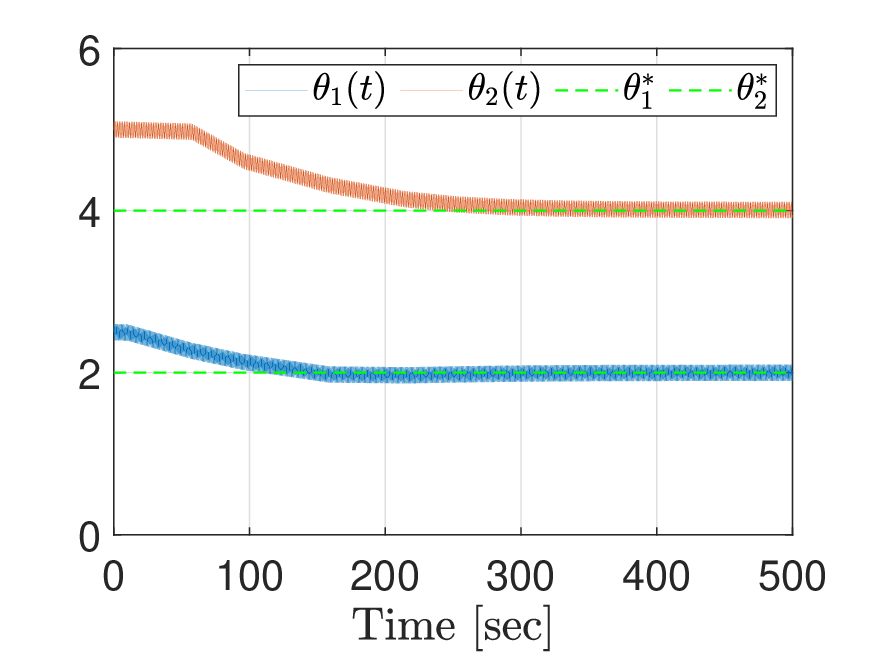}}    
	\\
	\subfigure[\underline{\bf ET-Gradient-ES:} control input, $u(t)$. \label{fig:ETGradientES_u}]{\includegraphics[width=4.1cm]{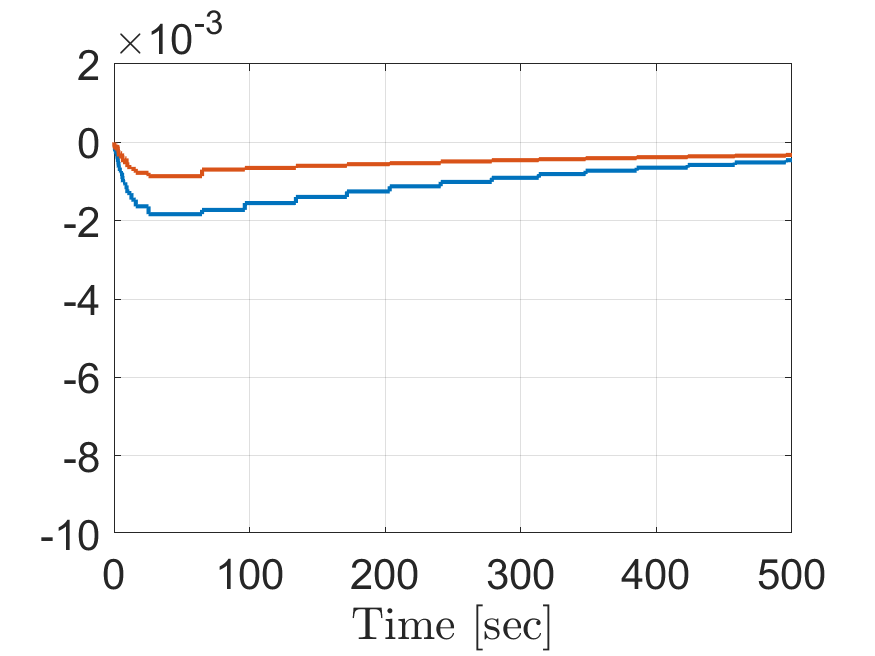}}
	~
	\subfigure[\underline{\bf ET-Newton-ES:} control input, $u(t)$. \label{fig::ETNewtonES_u}]{\includegraphics[width=4.1cm]{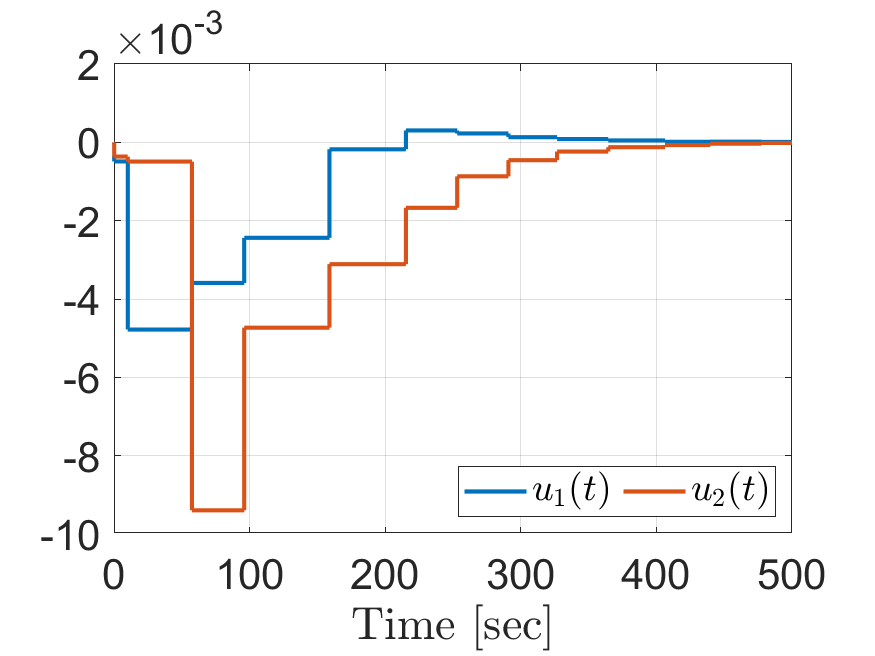}}
    \begin{picture}(300,0)(-47.5,-48)
    \includegraphics[width=2cm]{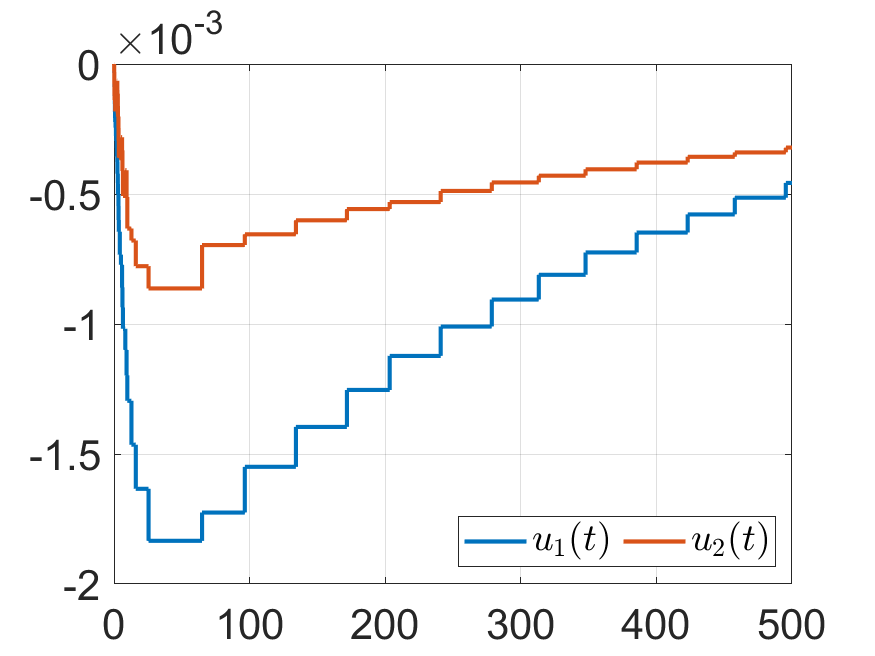}
    \end{picture}
    \begin{picture}(300,0)(-47.25,-205)
    \includegraphics[width=0.9cm]{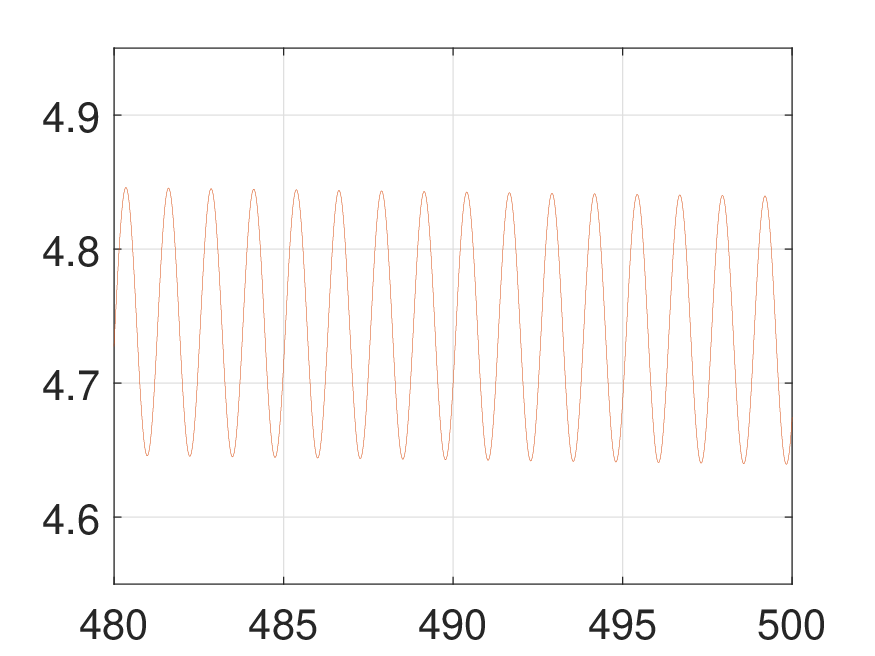}
    \end{picture}
    \begin{picture}(300,0)(-47.25,-190)
    \includegraphics[width=0.9cm]{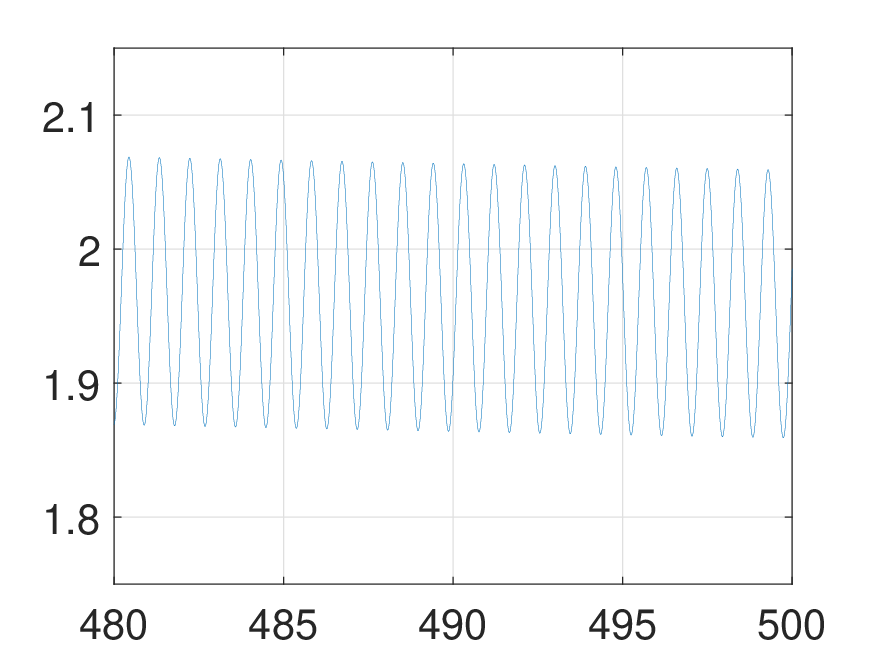}
    \end{picture}
    \begin{picture}(300,0)(-170.25,-223)
    \includegraphics[width=0.9cm]{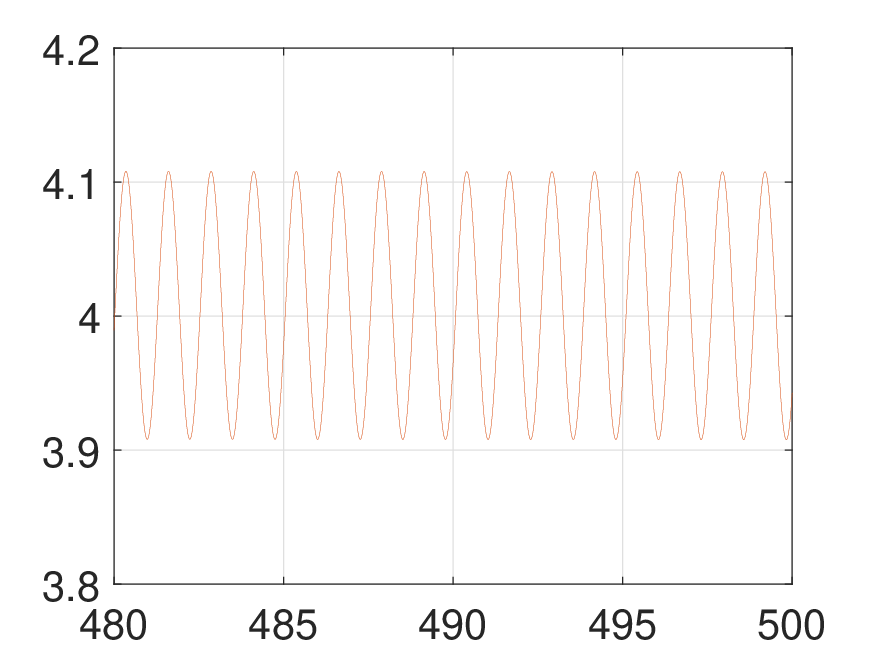}
    \end{picture}
    \begin{picture}(300,0)(-170.25,-212)
    \includegraphics[width=0.9cm]{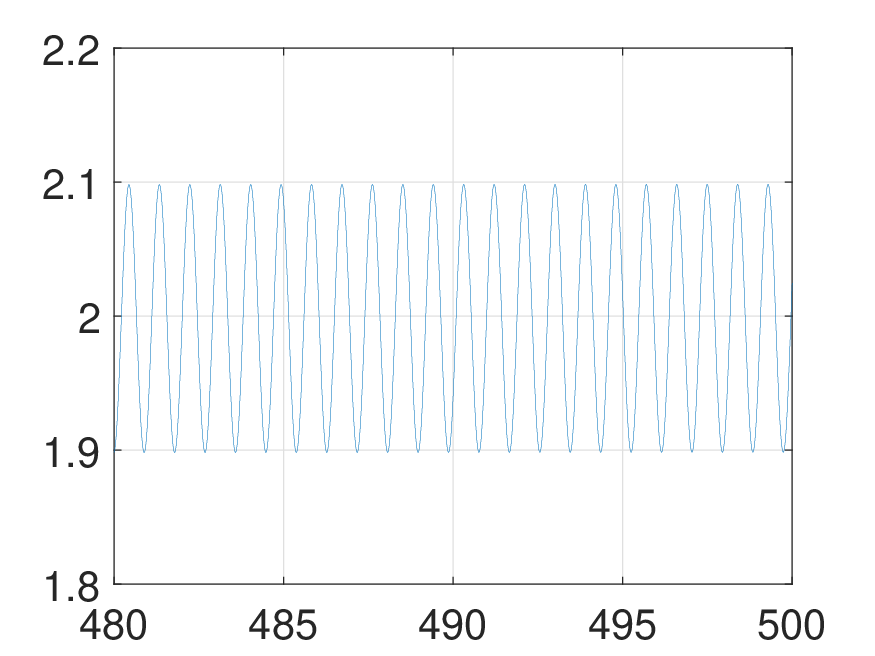}
    \end{picture}    

	\subfigure[\underline{\bf ET-Gradient-ES:} output of the map, $y(t)$. \label{fig:ETGradientES_y}]{\includegraphics[width=4.1cm]{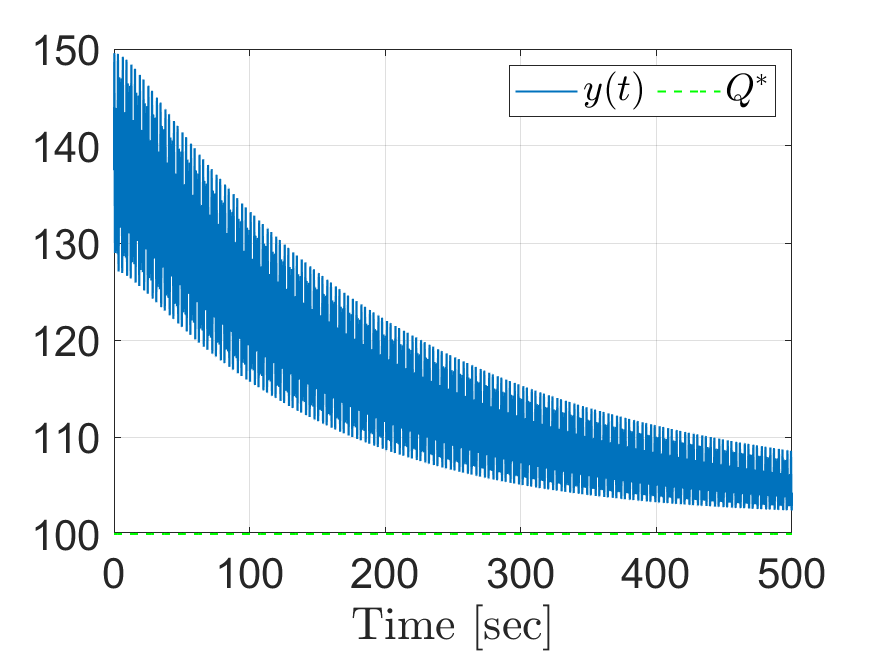}}
	~
	\subfigure[\underline{\bf ET-Newton-ES:} output of the map, $y(t)$. \label{fig:ETNewtonES_y}]{\includegraphics[width=4.1cm]{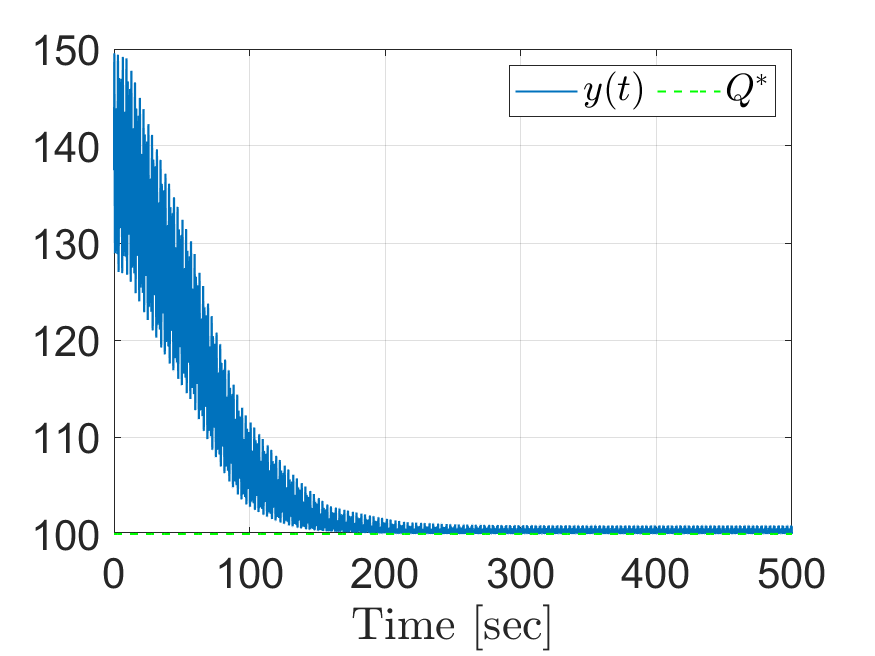}}
    \\
    \subfigure[\underline{\bf Gradient vs Newton:} Time evolution of the event-triggered controllers updates. \label{fig:updates}]{\includegraphics[width=4.1cm]{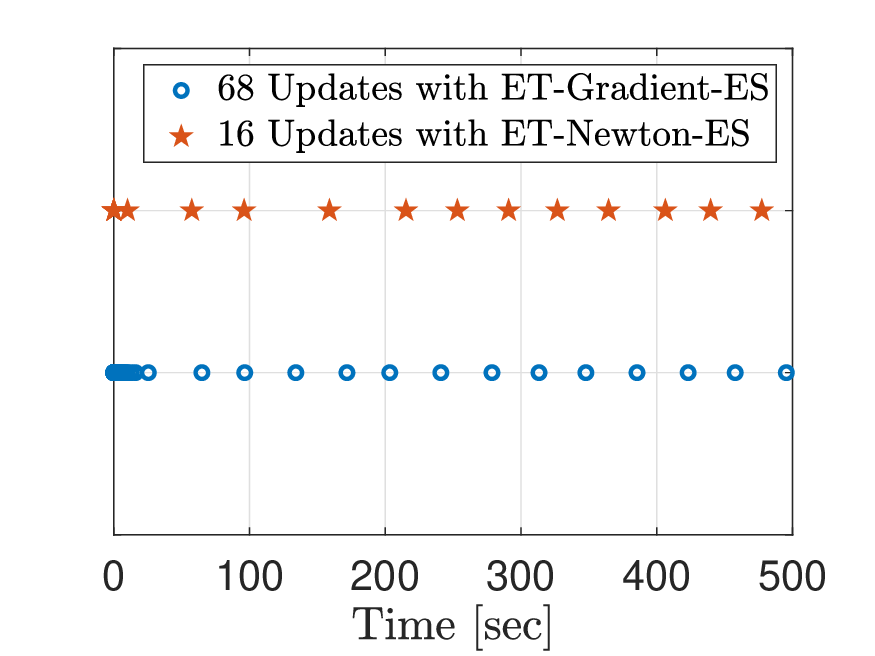}}
	~
    \subfigure[\underline{\bf ET-Newton-ES:} Hessian inverse estimate, $\Gamma(t)$. \label{fig:ETNewtonES_Gamma}]{\includegraphics[width=4.1cm]{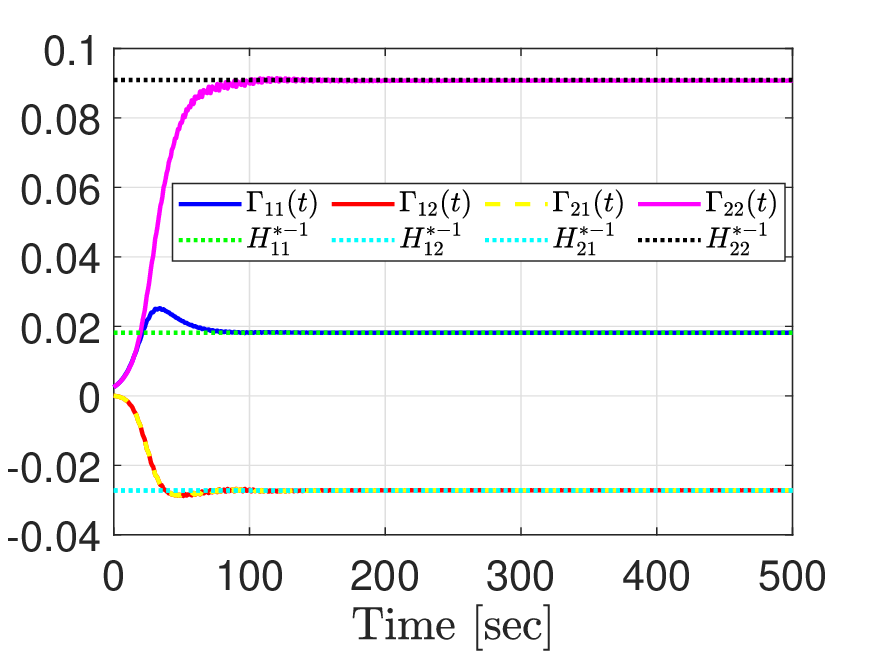}}
	\caption{Static and Dynamic Event-triggered Multivariable Extremum Seeking Systems. \label{fig:SD_ETESC}}
    \begin{picture}(300,0)(-54,-238)
    \includegraphics[width=1.75cm]{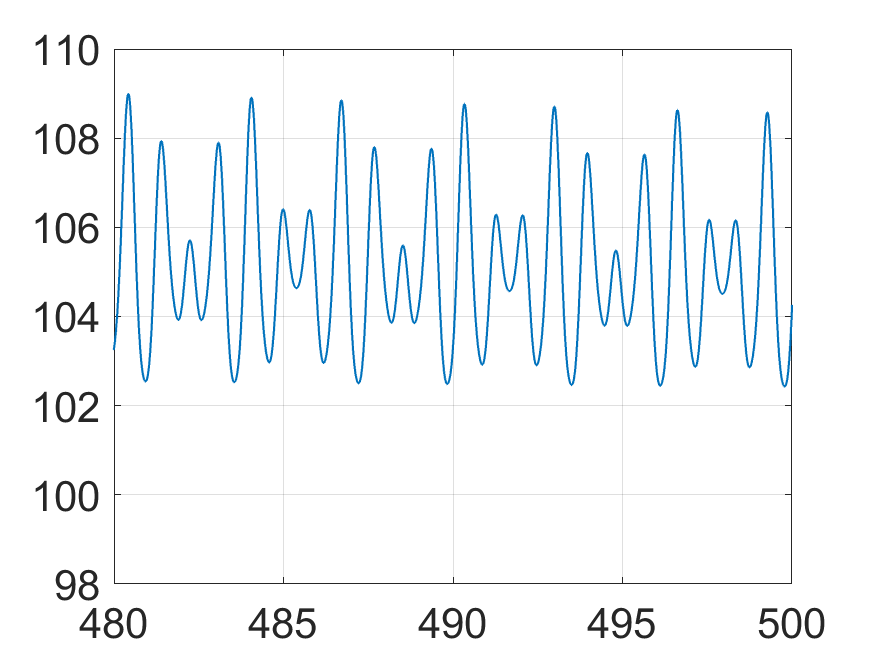}
    \end{picture}
    \begin{picture}(300,0)(-175.5,-240)
    \includegraphics[width=1.75cm]{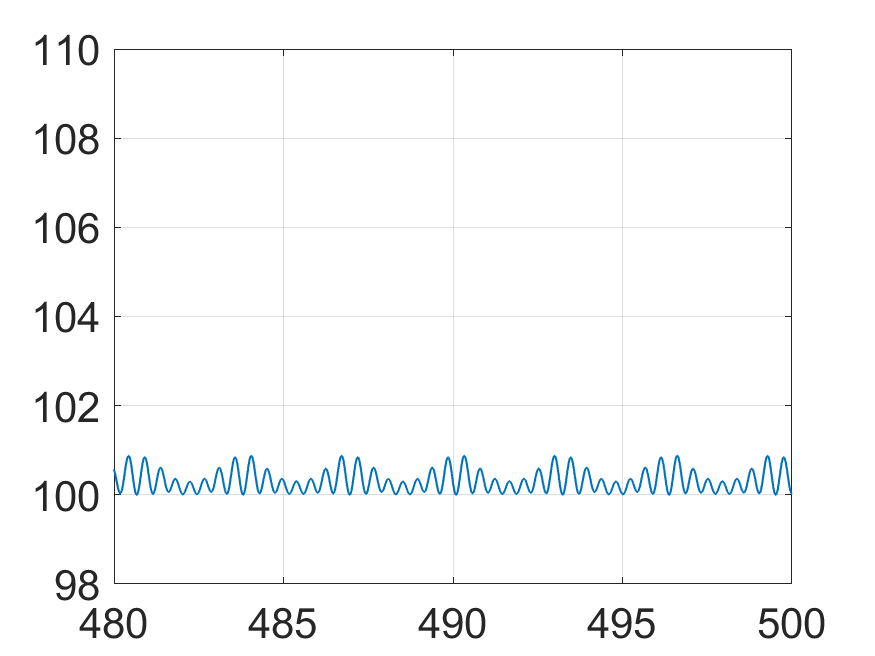}
    \end{picture}   
\end{figure}

The comparisons between the event-triggered gradient-based and Newton-based extremum seeking strategies highlights key performance distinctions, as illustrated in Figs.~\ref{fig:ETGradientES_theta} --\ref{fig:ETNewtonES_Gamma}. The evolution of the input variable \(\theta(t)\), displayed in Figs.~\ref{fig:ETGradientES_theta} and \ref{fig:ETNewtonES_theta}, further confirms the faster and more direct convergence of the Newton-based strategy, enabled by the use of second-order information curvature of the map. The control input \(u(t)\) generated by the gradient-based approach, shown in Fig.~\ref{fig:ETGradientES_u}, exhibits slower transients compared to the Newton-based method in Fig.~\ref{fig::ETNewtonES_u}, which produces sharper and more efficient control actions. In terms of output performance, both schemes drive the nonlinear output \(y(t)\) to its optimal value \(Q^\ast\), but Fig.~\ref{fig:ETNewtonES_y} shows a remarkable faster convergence compared to Fig.~\ref{fig:ETGradientES_y}. Notably, Fig.~\ref{fig:updates} reveals that the Newton-based method required only \textbf{16 control updates}, whereas the gradient-based strategy needed \textbf{68}, highlighting a significant reduction in actuation frequency. Finally, Fig.~\ref{fig:ETNewtonES_Gamma} illustrates the dynamic behavior of the inverse Hessian estimate \(\Gamma(t)\) in the Newton-based controller, which settles down after a brief transient phase, reinforcing the controller’s capability to adapt to the curvature of the cost function.

\section{Conclusions}

This paper introduced a static event-triggered control strategy tailored for multivariable Newton-based extremum seeking, with the goal of reducing the frequency of control signal updates while preserving fast and accurate convergence to the extremum of an unknown static multi-input nonlinear map. In contrast to traditional gradient-based extremum seeking schemes—whose convergence speed is inherently limited by the unknown and possibly ill-conditioned Hessian—the proposed strategy relies on a dynamic Riccati-based estimator of the inverse Hessian, thereby enabling user-defined convergence rates and decoupling performance from the curvature of the cost function.

By incorporating a static event-triggering condition into the Newton-based framework, the proposed approach ensures that control updates are executed only when necessary, significantly reducing actuation effort without sacrificing closed-loop stability or convergence accuracy. The theoretical stability guaranties are rigorously derived using averaging theory for systems with discontinuous right-hand sides, ensuring local exponential practical stability. Numerical simulations validate the analytical results and illustrate that the proposed event-triggered multivariable Newton-based strategy not only achieves faster convergence compared to its gradient-based counterpart, but also can requires significantly fewer control updates. This substantial reduction in actuation effort illustrates the efficiency of incorporating second-order information of the curvature map into the event-triggered framework, enabling more effective use of computational and communication resources.

The presented methodology advances the state-of-the-art in extremum seeking by bridging Newton-based optimization with sparse, feedback-dependent control policies, thus offering a viable and resource-efficient solution for real-time implementation in embedded or networked systems. Future work may explore dynamic event-triggered conditions, robustness to measurement noise, and extensions to time-varying or distributed optimization problems, including Nash equilibrium seeking in multi-agent settings.

\begin{ack}                               
The first and second authors thank the Brazilian funding agencies CAPES, CNPq and FAPERJ for the financial support.  
\end{ack}

\appendix

{\small
\section*{Appendix}

\section{Average for Systems with Discontinuous Right-Hand Sides}
\label{appendix_plotnikov}

It is important to clarify that
the following averaging result for differential inclusions with discontinuous right-hand sides by Plotnikov \cite{P:1979} takes into account discontinuities not in the periodic perturbations, but in the states. In our case, perturbations of the classical ESC scheme remains periodic while the discontinuities induced by the increasing sequence of event times  are not periodic. 

Instead of modeling \textcolor{black}{$e(\bar{t})$} as an additional state with its own dynamics or reset law (hybrid system approach), we treat \textcolor{black}{$e(\bar{t})$} as a discontinuous control input, constrained by the triggering condition (\ref{eq:tk+1_event}) such that its admissible values are restricted by the trigging inequality. From the triggering condition, we can define 
\begin{align}
    \Xi(\hat{G},\Gamma):=\sigma\|\Gamma(\bar{t})\hat{G}(\bar{t})\| \mathbb{-}\alpha\|\Gamma(\bar{t}_{\kappa})\hat{G}(\bar{t}_{\kappa})-\Gamma(\bar{t})\hat{G}(\bar{t})\|\,,
\end{align}
such that
\begin{align}
    \textcolor{black}{e(\bar{t})}=\begin{cases}
        ~~~~~~~~~~~~~~~0~~~~~~~~~~~~~\,, \quad  \Xi(\hat{G},\Gamma) < 0 \\
        \Gamma(\bar{t}_{\kappa})\hat{G}(\bar{t}_{\kappa})-\Gamma(\bar{t})\hat{G}(\bar{t})\,, \quad \Xi(\hat{G},\Gamma) \geq 0
    \end{cases}\,. \label{eq:eDI}
\end{align}
Notice, from (\ref{eq:tk+1_event}) and (\ref{eq:eDI}), between events,
\begin{align}
    \|\Gamma(\bar{t}_{\kappa})\hat{G}(\bar{t}_{\kappa})-\Gamma(\bar{t})\hat{G}(\bar{t})\| \leq \frac{\sigma}{\alpha}\|\Gamma(\bar{t})\hat{G}(\bar{t})\|\,.
\end{align}
Nevertheless, from the point of view of the discontinuous input (\ref{eq:eDI}), {\it i.e.}, even on the update instants, we have, for it component of \textcolor{black}{$e(\bar{t})$},
\begin{align}
    |e_{i}(\bar{t})| &\leq \|\textcolor{black}{e(\bar{t})}\| \nonumber \\
     &\leq \|\Gamma(\bar{t}_{\kappa})\hat{G}(\bar{t}_{\kappa})-\Gamma(\bar{t})\hat{G}(\bar{t})\| \nonumber \\
      &\leq  \frac{\sigma}{\alpha}\|\Gamma(\bar{t})\hat{G}(\bar{t})\|\,,
\end{align}
or, equivalently, 
\begin{align}
    -\frac{\sigma}{\alpha}\|\Gamma(\bar{t})\hat{G}(\bar{t})\|\leq e_{i}(\bar{t}) \leq \frac{\sigma}{\alpha}\|\Gamma(\bar{t})\hat{G}(\bar{t})\|\,, \quad \forall i \in \{1,\ldots,n\}\,.
\end{align}
Hence, once \textcolor{black}{$e(\bar{t})$} is viewed as an input, for all $i \in \{1,\ldots,n\}$, its \textcolor{black}{admissible set} is 
\begin{align}
    \textcolor{black}{\mathcal{E}(\Gamma,\hat{G})}\mathbb{:=}\left\{\textcolor{black}{e(\bar{t})}\in \mathbb{R}^{n}\Big|\mathbb{-}\frac{\sigma}{\alpha}\|\Gamma(\bar{t})\hat{G}(\bar{t})\|\mathbb{\leq} e_{i}(\bar{t}) \mathbb{\leq} \frac{\sigma}{\alpha}\|\Gamma(\bar{t})\hat{G}(\bar{t})\|\right\}\!. \label{eqE_DI}
\end{align}
The signal \textcolor{black}{$e(\bar{t})$} is allowed to be discontinuous, since it represents in some sense the sampling generated by the event-triggered mechanism whenever a new update occurs. However, despite these discontinuities, \textcolor{black}{$e(\bar{t})$} is constrained pointwise by the triggering condition such that the value of \textcolor{black}{$e(\bar{t})$} must belong to the admissible set (\ref{eqE_DI}). This pointwise constraint is the key property that allows the sampling error to be interpreted as a discontinuous input and absorbed into a differential inclusion.

With this interpretation, since $\Gamma(\bar{t})=\tilde{\Gamma}(\bar{t})+H^{\ast-1}$, the closed-loop becomes
\begin{align}
\frac{d\hat{G}}{d\bar{t}}(\bar{t})& \in \frac{1}{\omega} \mathcal{F}_{1}\left(\bar{t},\hat{G},\tilde{\theta},\tilde{\Gamma}\right) \,, \label{eq:dhatGdt_20260115_1}\\
\frac{d\tilde{\theta}}{d\bar{t}}(\bar{t})&\in \frac{1}{\omega} \mathcal{F}_{2}\left(\bar{t},\hat{G},\tilde{\theta},\tilde{\Gamma}\right)\,,\label{eq:dtildeThetadt_20260115_1} \\
\frac{d\tilde{\Gamma}}{d\bar{t}}(\bar{t})&\in \frac{1}{\omega} \mathcal{F}_{3}\left(\bar{t},\hat{G},\tilde{\theta},\tilde{\Gamma}\right)\,,\label{eq:dotTildeGamma_20260115_1}
\end{align}
where
\begin{align}
  &\mathcal{F}_{1}\left(\bar{t},\hat{G},\tilde{\theta},\tilde{\Gamma}\right) \mathbb{=} \mathbb{-}\left(H^{\ast}KH^{\ast-1}\mathbb{+}H^{\ast}K\tilde{\Gamma}(\bar{t}) \right. \nonumber \\
  &\left.\mathbb{+}\textcolor{black}{\Delta \! \mbox{\calligra H}^{~~~\ast}\!(\bar{t})}KH^{\ast-1}\mathbb{+}\textcolor{black}{\Delta \! \mbox{\calligra H}^{~~~\ast}\!(\bar{t})}K\tilde{\Gamma}(\bar{t}) \right)\hat{G}(\bar{t}) \nonumber \\
&\mathbb{+}\omega\textcolor{black}{\frac{d\Delta \! \mbox{\calligra H}^{~~~\ast}}{d\bar{t}}\!(\bar{t})}\tilde{\theta}(\bar{t})\mathbb{-}\left(\!H^{\ast}\mathbb{+}\textcolor{black}{\Delta \! \mbox{\calligra H}^{~~~\ast}\!(\bar{t})}\!\right)K\textcolor{black}{\mathcal{E}(\tilde{\Gamma},\hat{G})}\mathbb{+}\omega\textcolor{black}{\frac{d \Delta}{d\bar{t}}(\bar{t})}, \label{eq:dhatGdt_20260108_1} \\
&\mathcal{F}_{2}\left(\bar{t},\hat{G},\tilde{\theta},\tilde{\Gamma}\right) \mathbb{=}\mathbb{-}K\left(\! I_{n}\mathbb{+}\tilde{\Gamma}(\bar{t})H^{\ast}\mathbb{+}H^{\ast-1}\textcolor{black}{\Delta \! \mbox{\calligra H}^{~~~\ast}\!(\bar{t})} \right. \nonumber \\
&\left. \mathbb{+}\tilde{\Gamma}(\bar{t})\textcolor{black}{\Delta \! \mbox{\calligra H}^{~~~\ast}\!(\bar{t})}\!\right)\tilde{\theta}(\bar{t})\mathbb{-}K\textcolor{black}{\mathcal{E}(\tilde{\Gamma},\hat{G})}\mathbb{-}KH^{\ast\mathbb{-}1}\textcolor{black}{\Delta(\bar{t})}, \label{eq:dtildeThetadt_20260108_1} \\
&\mathcal{F}_{3}\left(\bar{t},\hat{G},\tilde{\theta},\tilde{\Gamma}\right) =\!\mathbb{-}\omega_{r}\tilde{\Gamma}(\bar{t})\mathbb{+}\omega_{r}\Upsilon(\bar{t},\tilde{\theta},\tilde{\Gamma}) \,.\label{eq:dotTildeGamma_20260108_1}
\end{align}

The closed-loop system \eqref{eq:dhatGdt_20250206_3}--\eqref{eq:dotTildeGamma_20250603_1} is governed by differential equations with discontinuous right-hand sides due to the triggering variable \textcolor{black}{$e(\bar{t})$}, and it is also periodic in $\bar{t}$. Moreover, by using the admissible set (\ref{eq:dhatGdt_20260108_1})--(\ref{eq:dotTildeGamma_20260108_1}), a solution of the non-autonomous system with differential equations (\ref{eq:dhatGdt_20250206_3})--(\ref{eq:dotTildeGamma_20250603_1}) is a solution of the non-autonomous system given by the differential inclusions (\ref{eq:dhatGdt_20260115_1})--(\ref{eq:dotTildeGamma_20260115_1}) also periodic in $\bar{t}$.

The differential inclusions (\ref{eq:dhatGdt_20260115_1})--(\ref{eq:dotTildeGamma_20260115_1}) are characterized by a small parameter $1/\omega$ as well as a time-periodic admissible set (\ref{eq:dhatGdt_20260115_1})--(\ref{eq:dotTildeGamma_20260115_1}) in $\bar{t}$ and, thereby, defining $x:=\begin{bmatrix} \hat{G}^{\top}\,, \tilde{\theta}^{\top}\,, \tilde{\Gamma}^{\top}\end{bmatrix}^{\top}$ and $\mathcal{F}:=\begin{bmatrix}  \mathcal{F}_{1}^{\top}\,,  \mathcal{F}_{2}^{\top}\,,  \mathcal{F}_{3}^{\top}\end{bmatrix}^{\top}$, we have the differential inclusion
\begin{align}
\dfrac{dx}{d\bar{t}}(\bar{t})&\in \dfrac{1}{\omega}\mathcal{F}\left(\bar{t},x\right)\,. \label{eq:dotX_event_DI}
\end{align} 
By following this procedure, we obtain a differential inclusion model of our event-triggered Newton-based extremum seeking scheme in the state variable. This model captures the family of trajectories consistent with the event-triggered implementation and provides a suitable representation for the subsequent analysis via the averaging method for differential inclusions proposed by Plotnikov \cite{P:1979}. Modeling event-triggered closed-loop dynamics as a differential inclusion has also been adopted in the literature, for instance by Delimpaltadakis and Mazo \cite{DM:2020,DM:2022}.

From \cite{P:1979}, let us consider the differential inclusion
\begin{align}
\frac{d x}{dt} \in \varepsilon \mathcal{F}(t,x)\,, \quad x(0)=x_{0}\,, \label{eq:A1}
\end{align}
where $x$ is an n-dimensional vector, $t$ is time, $\varepsilon$ is a small parameter, and $\mathcal{F}(t,x)$ is a multivalued function that is $T$-periodic in $t$ and puts in correspondence with each point $(t,x)$ of a certain domain of the ($n+1$)-dimensional space a compact set $\mathcal{F}(t,x)$ of the $n$-dimensional space. 
Let us put in correspondence with the inclusion  (\ref{eq:A1}) the average inclusion
\begin{align}
\frac{d \xi}{dt} \in \varepsilon \bar{\mathcal{F}}(\xi)\,, \quad \xi(0)=x_{0}\,, \label{eq:A2}
\end{align}
where $\bar{\mathcal{F}}(\xi)=\frac{1}{T}\int_{0}^{T}\mathcal{F}(\tau,\xi)d\tau$. 

\begin{theorem} \label{thm:A1}
Let a multivalued mapping $\mathcal{F}(t,x)$ be defined in the domain $Q\left\{t\geq 0\,, x\in D\subset\mathbb{R}^{n}\right\}$ and let in this domain the set $\mathcal{F}(t,x)$ be a nonempty compactum for all admissible values of the arguments and the mapping $\mathcal{F}(t,x)$ be continuous and uniformly bounded and satisfy the Lipschitz condition with respect to $x$ with a constant $\lambda$, {\it i.e.}, $\mathcal{F}(t,x) \subset S_{M}(0)$, $\delta(\mathcal{F}(t,x')-\mathcal{F}(t,x''))\leq \lambda \|x'-x''\|$, where $\delta(P,Q)$ is the Hausdorff distance between the sets $P$ and $Q$, {\it i.e.}, $\delta(P,Q)=\min\left\{d|P\subset S_{d}(Q), Q \subset S_{d}(P)\right\}$, $S_{d}(N)$ being the d-neighborhood of a set $N$ in the space $\mathbb{R}^{n}$; the mapping $\mathcal{F}(t,x)$ be $T$-periodic in $t$; for all $x_{0}\in D'\subset D$ the solutions of inclusion (\ref{eq:A2}) lie in  the domain $D$ together with a certain $\rho$-neighborhood. Them for each $L>0$ there exist $\varepsilon^{0}(L)>0$ and $c(L)>0$ such that for $\varepsilon \in (0\,,\varepsilon^{0}\rbrack$ and $t \in \lbrack 0, L\varepsilon^{-1}\rbrack$:
\begin{enumerate}
\item for each solution $x(t)$ of the inclusion (\ref{eq:A1}) there exists a solution $\xi(t)$ of the inclusion (\ref{eq:A2}) such that 
\begin{align}
\|x(t)-\xi(t)\|\leq c\varepsilon =\mathcal{O}(\varepsilon); \label{eq:A4}
\end{align}
\item for each solution $\xi(t)$ of the inclusion (\ref{eq:A2}) there exists a solution $x(t)$ of the inclusion (\ref{eq:A1}) such that  the inequality (\ref{eq:A4}) holds.
\end{enumerate}
Thus, the following estimate is valid: $\delta(\bar{R}(t),R'(t))\leq c\varepsilon =\mathcal{O}(\varepsilon)$, 
where $\bar{R}(t)$ is a section of the family of solutions of the inclusion (\ref{eq:A2}) and $R'(t)$ is the closure of the section $R(t)$ of the family of solutions of the inclusion (\ref{eq:A1}).
\end{theorem}

\begin{theorem} \label{thm:A2}
Let all the conditions of Theorem~\ref{thm:A1} and also the following condition be fulfilled: the $R$-solution $\bar{R}(t)$ of inclusion (\ref{eq:A2}) is uniformly asymptotically stable. Then there exist $\varepsilon^{0}>0$ and $c>0$ such that for $0< \varepsilon \leq \varepsilon^{0}$
\begin{align}
\delta(\bar{R}(t),R'(t))\leq c\varepsilon = \mathcal{O}(\varepsilon)\,, \quad \forall t\geq 0. \label{eq:A6}
\end{align}
\end{theorem}
}

\begin{thebibliography}{10}
\expandafter\ifx\csname url\endcsname\relax
  \def\url#1{\texttt{#1}}\fi
\expandafter\ifx\csname urlprefix\endcsname\relax\def\urlprefix{URL }\fi
\expandafter\ifx\csname href\endcsname\relax
  \def\href#1#2{#2} \def\path#1{#1}\fi

\bibitem{A:1957}
T.~Apostol.
\newblock {\em Mathematical Analysis - A Modern Approach to Advanced Calculus}.
\newblock Addison-Wesley Publishing Company, Massachusetts, 1957.

\bibitem{AK:2003}
K.~B.~Ariyur and M.~Krsti{\' c}.
\newblock {\em Real-Time Optimization by Extremum-Seeking Control}.
\newblock Wiley, Canada, 2003.

\bibitem{Atta2015}
K.~T.~Atta, A.~Johansson, and T.~Gustafsson.
\newblock Extremum seeking control based on phasor estimation.
\newblock {\em Systems and Control Letters}, 85:37--45, 2015.

\bibitem{DM:2020}
G.~Delimpaltadakis and M.~Mazo, Jr.
\newblock Region-Based Self-Triggered Control for Perturbed and Uncertain Nonlinear Systems.
\newblock {\em IEEE Transactions on Automatic Control}, 2(8):757--768, 2021.

\bibitem{DM:2022}
G.~Delimpaltadakis and M.~Mazo, Jr.
\newblock Traffic Abstractions of Nonlinear Event-Triggered Control Systems with Disturbances and Uncertainties.
\newblock {\em arXiv preprint}, \url{https://www.arxiv.org/pdf/2010.12341v1}, 2022.

\bibitem{Espitia2020}
N.~Espitia.
\newblock Observer-based event-triggered boundary control of a linear $2\times 2$ hyperbolic systems.
\newblock {\em Systems and Control Letters}, 138, 2020.

\bibitem{FKB:2012}
P.~Frihauf, M.~Krsti{\' c}, and T.~Ba{\c s}ar.
\newblock Nash equilibrium seeking in noncooperative games.
\newblock {\em IEEE Transactions on Automatic Control}, 57:1192--1207, 2012.

\bibitem{GKN:2012}
A.~Ghaffari, M.~Krsti{\'c}, and D.~Ne{\u s}ic.
\newblock Multivariable {Newton}-based extremum seeking.
\newblock {\em Automatica}, 48:1759--1767, 2012.

\bibitem{G:2014}
A.~Girard.
\newblock Dynamic triggering mechanism for event-triggered control.
\newblock {\em IEEE Transactions on Automatic Control}, 60:1992--1997, 2014.

\bibitem{HJT:2012}
W.~P.~M.~H.~Heemels, K.~H.~Johansson, and P.~Tabuada.
\newblock An introduction to event-triggered and self-triggered control.
\newblock In {\em 2012 IEEE 51st IEEE Conference on Decision and Control (CDC)}, 3270--3285, 2012.

\bibitem{HNX:2007}
J.~P.~Hespanha, P.~Naghshtabrizi, and Y.~Xu.
\newblock A survey of recent results in networked control systems.
\newblock {\em Proceedings of the IEEE}, 95(1):138--162, 2007.

\bibitem{K:2002}
H.~K.~Khalil.
\newblock {\em Nonlinear Systems}.
\newblock Prentice Hall, Upper Saddle River, New Jersey, 2002.

\bibitem{Krstic2000}
M.~Krstić.
\newblock Performance improvement and limitations in extremum seeking control.
\newblock {\em Systems and Control Letters}, 39(5):313--326, 2000.

\bibitem{K:2014}
M.~Krsti{\'c}.
\newblock Extremum seeking control.
\newblock In J.~Baillieul and T.~Samad (eds.), {\em Encyclopedia of Systems and Control}, 1:1--5. Springer, London, 1st ed., 2014.

\bibitem{KW:2000}
M.~Krsti{\' c} and H.-H.~Wang.
\newblock Stability of extremum seeking feedback for general nonlinear dynamic systems.
\newblock {\em Automatica}, 36:595--601, 2000.

\bibitem{L:1922}
M.~Leblanc.
\newblock Sur l'{\' e}lectrification des chemins de fer au moyen de courants alternatifs de fr{\' e}quence {\' e}lev{\' e}e.
\newblock {\em Revue G{\' e}n{\' e}rale de l'Electricit{\' e}}, XII(8):275--277, 1922.

\bibitem{Mazenc2022}
F.~Mazenc, M.~Malisoff, and C.~Barbalata.
\newblock Event-triggered control for continuous-time linear systems with a delay in the input.
\newblock {\em Systems and Control Letters}, 159, 2022.

\bibitem{MMB:2010}
W.~H.~Moase, C.~Manzie, and M.~J.~Brear.
\newblock Newton-like extremum-seeking for the control of thermoacoustic instability.
\newblock {\em IEEE Transactions on Automatic Control}, 55(9):2094--2105, 2010.

\bibitem{OLIVEIRA2015304}
T.~R.~Oliveira and M.~Krsti{\' c}.
\newblock Newton-based extremum seeking under actuator and sensor delays.
\newblock {\em IFAC-PapersOnLine}, 48(12):304--309, 2015. (12th IFAC Workshop on Time Delay Systems, TDS 2015).

\bibitem{TRoux:2022}
T.~R.~Oliveira and M.~Krsti{\' c}.
\newblock {\em Extremum Seeking through Delays and PDEs}.
\newblock Society for Industrial and Applied Mathematics, 2022.

\bibitem{TRoux:2021dec}
T.~R.~Oliveira, V.~H.~P.~Rodrigues, M.~Krsti{\' c}, and T.~Ba{\c s}ar.
\newblock Nash equilibrium seeking in heterogeneous noncooperative games with players acting through heat {PDE} dynamics and delays.
\newblock In {\em 2021 60th IEEE Conference on Decision and Control (CDC)}, 1167--1173, 2021.

\bibitem{TRoux:2021may}
T.~R.~Oliveira, V.~H.~P.~Rodrigues, M.~Krsti{\' c}, and T.~Ba{\c s}ar.
\newblock Nash equilibrium seeking with players acting through heat {PDE} dynamics.
\newblock In {\em 2021 American Control Conference (ACC)}, 684--689, 2021.

\bibitem{ORKB:2021}
T.~R.~Oliveira, V.~H.~P.~Rodrigues, M.~Krsti{\' c}, and T.~Ba{\c s}ar.
\newblock Nash equilibrium seeking with players acting through heat {PDE} dynamics.
\newblock {\em Journal of Optimization Theory and Applications}, 191:700--735, 2021.

\bibitem{TRO:2020}
P.~Paz, T.~R.~Oliveira, A.~V.~Pino, and A.~P.~Fontana.
\newblock Model-free neuromuscular electrical stimulation by stochastic extremum seeking.
\newblock {\em IEEE Transactions on Control Systems Technology}, 28:238--253, 2020.

\bibitem{P:1979}
V.~A.~Plotnikov.
\newblock Averaging of differential inclusions.
\newblock {\em Ukrainian Mathematical Journal}, 31:454--457, 1980.

\bibitem{VHPR:2022}
V.~H.~P.~Rodrigues, L.~Hsu, T.~R.~Oliveira, and M.~Diagne.
\newblock Event-triggered extremum seeking control.
\newblock {\em IFAC-PapersOnLine}, 55(12):555--560, 2022.

\bibitem{VHPR:2023b}
V.~H.~P.~Rodrigues, L.~Hsu, T.~R.~Oliveira, and M.~Diagne.
\newblock Dynamic event-triggered extremum seeking feedback.
\newblock {\em IFAC-PapersOnLine}, 56(2):10307--10314, 2023.

\bibitem{VHPR:2023a}
V.~H.~P.~Rodrigues, T.~R.~Oliveira, L.~Hsu, M.~Diagne, and M.~Krstic.
\newblock Event-Triggered and Periodic Event-Triggered Extremum Seeking Control.
\newblock {\emph Automatica}, 174:112161, 2025.

\bibitem{TRO:2019}
T.~Roux{-}Oliveira, L.~R.~Costa, A.~V.~Pino, and P.~Paz.
\newblock Extremum seeking-based adaptive {PID} control applied to neuromuscular electrical stimulation.
\newblock {\em Annals of the Brazilian Academy of Sciences}, 21:1--20, 2019.

\bibitem{Scheinker2014}
A.~Scheinker and M.~Krstić.
\newblock Extremum seeking with bounded update rates.
\newblock {\em Systems and Control Letters}, 63:25--31, 2014.

\bibitem{SD:2009}
M.~Stankovi{\'c} and D.~Stipanovi{\'c}.
\newblock Discrete time extremum seeking by autonomous vehicles in a stochastic environment.
\newblock In {\em Proceedings of the 48th IEEE Conference on Decision and Control held jointly with 2009 28th Chinese Control Conference, CDC/CCC 2009}, 4541--4546, 2009.

\bibitem{T:2007}
P.~Tabuada.
\newblock Event-triggered real-time scheduling of stabilizing control tasks.
\newblock {\em IEEE Transactions on Automatic Control}, 52:1680--1685, 2007.

\bibitem{Wang2012}
X.~Wang, Y.~Sun, and N.~Hovakimyan.
\newblock Asynchronous task execution in networked control systems using decentralized event-triggering.
\newblock {\em Systems and Control Letters}, 61(9):936--944, 2012.

\bibitem{W:1971}
F.~Warner.
\newblock {\em Foundations of Differentiable Manifolds and Lie Groups}.
\newblock Scott Foresman and Company, Chicago, Illinois, 1971.

\bibitem{zhang2012extremum}
C.~Zhang and R.~Ord{\'o}{\~n}ez.
\newblock {\em Extremum-Seeking Control and Applications: A Numerical Optimization-Based Approach}.
\newblock Springer, New York, NY, 2012.

\bibitem{ZHGDDYP:2020}
X.-M.~Zhang, Q.-L.~Han, X.~Ge, D.~Ding, L.~Ding, D.~Yue, and C.~Peng.
\newblock Networked control systems: A survey of trends and techniques.
\newblock {\em IEEE/CAA Journal of Automatica Sinica}, 7(1):1--17, 2020.

\bibitem{Zhang2017}
P.~Zhang, T.~Liu, and Z.~P.~Jiang.
\newblock Input-to-state stabilization of nonlinear discrete-time systems with event-triggered controllers.
\newblock {\em Systems and Control Letters}, 103:16--22, 2017.

\bibitem{ZF:2022}
Y.~Zhu and E.~Fridman.
\newblock Extremum seeking via a time-delay approach to averaging.
\newblock {\em Automatica}, 135:109965, 2022.

\bibitem{ZFO:2022}
Y.~Zhu, E.~Fridman and T.~R.~Oliveira.
\newblock Sampled-data extremum seeking with constant delay: a time-delay approach.
\newblock {\em IEEE Transactions on Automatic Control}, 68(1), 432--439, 2022.

\end{thebibliography}
\end{document}